\documentclass[10pt, twoside, fleqno, a4paper]{article} 

\usepackage{stmaryrd}
\usepackage{amsmath}
\usepackage{a4}
\usepackage[francais, english]{babel}
\usepackage{euscript}
\usepackage{amsfonts}
\usepackage{amssymb}
\usepackage{xypic}
\usepackage[dvips]{graphics}
\usepackage{enumerate}
\usepackage{array}
\usepackage{theorem}
\usepackage{epic}
\usepackage{color}
\xyoption{all}
\usepackage[T1]{fontenc}
\usepackage{makeidx}
\usepackage{float}

\DeclareMathAlphabet{\mathpzc}{OT1}{pzc}{m}{it}

\def\binome#1#2{
\left(\scriptstyle{#1}\atop \scriptstyle{#2}\right )}

\newtheorem{theorem}{Th\'eor\`eme}[subsection]
\newtheorem{proposition}[theorem]{Proposition}
\newtheorem{corollary}[theorem]{Corollaire}
\newtheorem{lemma}[theorem]{Lemme}

\theorembodyfont{\upshape}
\newtheorem{definition}[theorem]{D\'efinition}
\newtheorem{remark}[theorem]{Remarque}
\newtheorem{exemple}[theorem]{Exemple}

\numberwithin{equation}{subsection}

\newcommand{\qed}{\hfill $\Box$ \medskip}

\newcommand{\NN}{{\mathbb N}}
\newcommand{\ZZ}{{\mathbb Z}}

\renewcommand{\O}{{\cal O}}

\newcommand{\C}{{\cal C}}

\newcommand{\Spec}{{\operatorname{Spec}\kern 1pt}}
\newcommand{\Spf}{{\operatorname{Spf}\kern 1pt}}

\def\buildrel#1\over#2{\mathrel{\mathop{\kern 0pt#2}\limits^#1}}

\newcommand{\Ext}{{\mathrm{Ext}}}
\newcommand{\ext}{{\mathpzc{Ext}}}
\newcommand{\Hom}{{\mathrm{Hom}}}
\renewcommand{\hom}{{\mathpzc{Hom}}}
\renewcommand{\H}{{\mathrm{H}}}
\newcommand{\R}{{\mathrm{R}}}
\renewcommand{\Im}{\operatorname{Im}\kern 1pt}
\newcommand{\coker}{\operatorname{coker}\kern 1pt}

\newcommand{\card}{{\mathrm{card}}}

\renewcommand{\epsilon}{\varepsilon}
\newcommand{\Aut}{{\mathrm{Aut}}}
\newcommand{\Frac}{\mathrm{Frac}}
\newcommand{\preuve}{\noindent{\it Preuve : }}

\newcommand{\Ar}{{\mathbf{Ar}\kern 0.5pt}}

\newcommand{\cd}{\operatorname{cd}\kern 1pt}

\renewcommand{\tilde}{\widetilde}
\renewcommand{\hat}{\widehat}
\newcommand{\Tr}{\textrm{Tr}}

\newcommand{\X}{{\mathcal X}}

\newcommand{\p}{{\mathfrak p}}
\newcommand{\q}{{\mathfrak q}}

\newcommand{\A}{{\mathfrak A}}

\newcommand{\m}{{\mathfrak m}}

\newcommand{\Hc}{{\mathcal H}}

\newcommand{\dif}{{\mathfrak d}}

\def\d#1{{\partial \over \partial #1}}

\def\der#1{{\partial_{#1}}}

\newcommand{\cart}{\ar@{}[dr] |{\square}}
\newcommand{\comm}{\ar@{}[dr] |{\circlearrowleft}}
\begin{document}

\newcommand{\Fg}{{\mathfrak F}}
\newcommand{\Gg}{{\mathfrak G}}
\newcommand{\Na}{{\mathfrak{N}}}

\title{Déformations équivariantes des courbes stables I  Étude
cohomologique}

\author{Sylvain Maugeais}

\maketitle

\selectlanguage{english}

\begin{abstract}
Let $k$ be a field, $C\to \Spec k$ be a stable curve and let $G$ be a finite group
acting faithfully on the curve $C\to \Spec k$. In this article, we compute the 
vector space $\Ext^1_G(\Omega_{C/k}, \O_C)$, the sheaf $\Omega_{C/k}$ being the sheaf
of relative differentials. This vector space is naturally isomorphic to the 
set of first order $G$-equivariant deformations of $C$. The computation we do here 
will be used in a future article.  
\end{abstract}

\selectlanguage{francais}

\section{Introduction}

Soient $k$ un corps algébriquement clos de caractéristique $p \ge 0$, $C\to
\Spec k$ une courbe stable et $G$ un groupe fini agissant fidèlement
sur $C$. Le but de cet article est de décrire l'ensemble des
déformations équivariantes du couple $(C, G)$. Il est bien connu
que cet ensemble est naturellement isomorphe à
$\Ext^{1}_{G} (\Omega_{C/k}, \O_C)$ où $\Omega_{C/k}$ désigne le
faisceau des différentielles relatives.

Lorsque $p$ ne divise pas l'ordre de $G$, il est assez facile de
décrire l'espace vectoriel $\Ext^{1} (\Omega_{C/k}, \O_C)$ (cf. par
exemple \cite{Tuffery}). Les
vrais problèmes apparaissent lorsque $p$ divise l'ordre de $G$. Sous
des hypothèses supplémentaires, ce problème a déjà été étudié  par
plusieurs auteurs. Parmi ceux-ci, citons Laudal et Lønsted dans
\cite{Laudal_Lonsted} qui ont ainsi prouvé la lissité de l'espace des
modules des courbes hyperelliptiques en caractéristiques $2$ ; Bertin
et Mézard dans \cite{BertinMezard} qui ont étudié le cas général des
courbes lisses munies d'une action d'un $p$-groupe cyclique ;
Cornelissen et Kato dans \cite{CornelissenKato} qui travaillent avec des
groupes généraux et supposent que les courbes sont lisses et
ordinaires.

Lorsque $C\to \Spec k$ est une courbe stable, l'étude se complique
pour plusieurs raisons. La première est, bien sûr, l'existence de
singularité et peut être traitée de manière purement locale. La
deuxième vient du fait que $G$ peut opérer trivialement sur certaines
composantes irréductibles. Ce dernier problème est donc plus subtil
car il est de nature globale.

Décrivons maintenant le plan de cet article. La deuxième section fixe
le cadre de l'étude et énonce des résultats classiques concernant les
déformations équivariantes. Dans la troisième section,  nous
nous concentrons sur le cas local. Tout d'abord nous précisons
certains résultats obtenus par Bertin et Mézard dans le cas lisse.
Plus précisément, dans le cas où $G$ est un groupe cyclique, nous 
décrivons l'existence d'une base de l'espace des champs de vecteurs ayant une trace 
nulle sous l'action de $G$ (cf. Théorème \ref{ElementNonNulTraceNulle}). 
Ce résultat est fondamental pour la suite 
et il apparaîtra à plusieurs endroits. Nous passons ensuite au cas
singulier. L'idée est de relier certaines déformations équivariantes à
l'épaississement de la singularité et les autres aux déformations
équivariantes de la normalisée de $C$ (cf. suite exacte \eqref{EquationPointDouble}).

Dans la deuxième partie, nous étudions le cas global. Notons
$\pi:C\to C/G$ le morphisme quotient. Nous étudions chaque terme du
diagramme suivant (dont les lignes et les colonnes sont exactes) :

\begin{footnotesize}
$$\xymatrix@C0.3cm{
 & & & 0 \ar[d] & \\
 & & & \H^0(C/G, \R^1 \pi_*^G \Omega_{C/k}^\vee)\ar[d] & \\
0 \ar[r] & \H^1(C/G, \ext^0_G(\Omega_{C/k}, \O_C)) \ar[r] & \Ext^1_G(\Omega_{C/k}, \O_C) \ar[r]\ar[rd]_{\textit{ép}} & \H^0(C/G, \ext^1_G(\Omega_{C/k}, \O_C)) \ar[r]\ar[d] & 0 \\
 & & & \H^0(C/G, \pi_*^G \ext^1(\Omega_{C/k}, \O_C)) & 
}$$
\end{footnotesize}

L'étude du terme de gauche et du terme du haut (cf. Théorème
\ref{LienNormaliseGlobale}) est faite en se 
ramenant aux déformations équivariantes de la normalisée et en
utilisant les résultats locaux obtenus dans la troisième
section. L'image de ${\textit{ép}}$ est le terme le plus compliqué, son étude 
nécessite des méthodes de recollement formel (contrairement à l'étude
des autres termes qui repose seulement sur la cohomologie
équivariante) et fait l'objet du théorème \ref{EpaississementPointdouble}. 

Cet article sera suivi d'une étude de l'espace des modules des courbes
stables munies d'une action d'un groupe. Cette étude utilisera les
résultats démontrés dans le présent article et les illustrera. 

\medskip
\noindent{\bf Notations :} Soient $X$ un schéma, $G$ un groupe
agissant sur $X$ et $\p$ un point de $X$ (non nécessairement fermé). 
On notera $D_\p\subseteq G$ le stabilisateur de $\p$, $T_\p$ le
noyau de l'homomorphisme $D_\p \to \Aut(\Spec \O_{X, \p})$ et
$G_\p:=D_\p/T_\p$. 
 
Les formes différentielles de $k[[x]]$ seront notées $f(x) dx$, les champs
de vecteurs $f (x)\d x$ et la dérivation par rapport à la variable $x$ sera
notée $\der x$. 

La cohomologie équivariante est le principal outil utilisé dans cet
article, les résultats utilisés pourront être trouvés dans \cite{Tohoku}.

\section{Rappels sur les déformations équivariantes}

\begin{definition}
Soient $G$ un groupe fini et $Y$ un schéma. Un \textit{$Y[G]$-schéma}
(ou $A[G]$-schéma si $Y=\Spec A$) est un couple $(X\to Y, i)$ où $X$
est un $Y$-schéma et $i:G \to \Aut_Y(X)$ est un homomorphisme de
groupes. S'il  n'y a pas de confusion possible, on notera aussi $(X\to
Y, G)$, voire $(X, G)$, le couple $(X\to Y, i)$. Les $Y[G]$-schémas
forment une catégorie (les morphismes de $Y[G]$-schémas étant les
morphismes de $Y$-schémas commutant à l'action de $G$). 
\end{definition}

Soient $k$ un corps, $G$ un groupe fini et $(X, G)$ un
$k[G]$-schéma. Notons $\A$ la catégorie des anneaux locaux artiniens
de corps résiduel $k$. Une \textit{déformation $G$-équivariante} de
$(X, G)$ est un diagramme cartésien (dans la catégorie des
$\ZZ[G]$-schémas) 
$$\xymatrix{
(X, G) \ar[r]\ar[d] & (\X, G) \ar[d]\\
\Spec \ k \ar[r] & \Spec A
}$$
o\`u $A$ est un objet de $\A$ et $\X \rightarrow \Spec A$ est un morphisme plat. Deux déformations équivariantes $(\X, G)$ et $(\X', G)$ au-dessus de $A$ sont dites isomorphes s'il existe un $A$-morphisme équivariant $f:(\X, G) \rightarrow (\X', G)$ qui induit l'identité sur $(X, G)$ (remarquons qu'alors $f$ induit un isomorphisme $\X \to \X'$ d'après \cite{Deformation}, Lemma 3.3).

\begin{definition}\label{DefinitionRelevement}
Soient $A'$ un objet de $\A$ d'idéal maximal $\m_{A'}$ et ${\mathfrak
a} \subset A'$ un idéal tel que ${\m_{A'} \mathfrak a}=0$ ; notons
$A:=A'/{\mathfrak a}$. Soit $(\X, G)$ une déformation $G$-équivariante
de $(X, G)$ au-dessus de $\Spec A$. Un \textit{relèvement équivariant}
de $(\X, G)$ à $A'$ est une déformation équivariante $(\X', G)$ de
$(X, G)$ au-dessus de $\Spec A'$ telle que $(\X', G)\times _{\Spec A'}
\Spec A \cong (\X, G)$. 
\end{definition}

On a alors un résultat de classification des relèvements équivariants
en terme de groupes de cohomologie. 

\begin{theorem}[Illusie]
\label{EspaceTangent}
Soient $k$ un corps, $G$ un groupe fini, $X$ un $k$-schéma noe\-thérien
réduit et localement d'intersection complète muni d'une action de
$G$. Notons $\Omega_{X/k}$ le faisceau des différentielles
relatives. Reprenons les notations de la définition
\ref{DefinitionRelevement} et supposons qu'il existe une déformation
$G$-équivariante $(\X, G)$ de $(X, G)$ au-dessus de $\Spec A$. 
\begin{enumerate}
\item Si $(\X', G)$ est un relèvement équivariant de $(\X, G)$ à
$\Spec A'$ alors le groupe des automorphismes (de relèvement
équivariant) de $(\X', G)$ est canoniquement isomorphe à $${\mathfrak
a} \otimes_k \Hom_{\O_X}(\Omega_{X/k}, \O_X)^G.$$ 
\item S'il existe un relèvement de $(\X, G)$ à $\Spec A'$ alors il y a
une action canonique de ${\mathfrak a} \otimes_k
\Ext^1_G(\Omega_{X/k}, \O_X)$ sur l'ensemble des classes d'équivalence
de relèvements (modulo isomorphisme) de $(\X, G)$ à $\Spec A'$,
faisant de ce dernier un espace principal homogène. 
\end{enumerate}
De plus, si $X$ est le spectre d'un anneau local complet alors
$\Omega_{X/k}$ peut être remplacé par le module des différentielles
complété. 
\end{theorem}

Ce résultat est bien connu, on peut en trouver une preuve dans
\cite{IllusieI, IllusieII} où dans \cite{DeformationWewers}.

L'ensemble des classes d'isomorphies de relèvements de $(X, G)$ à 
$k[\varepsilon]/(\varepsilon^{2})$ (que nous appellerons 
\textit{déformations du premier ordre}) s'identifie donc naturellement à 
$\Ext^1_G(\Omega_{X/k}, \O_X)$ car il existe une déformation
équivariante au-dessus de $\Spec k[\epsilon]/(\epsilon^2)$ (à savoir
la déformation triviale $(X, G) \times_{\Spec k} \Spec
k[\epsilon]/(\epsilon^2)$). 

Nous nous limiterons à l'étude de l'action des $p$-groupes en
caractéristique $p >0$, cela étant en partie justifié par le lemme suivant.
\begin{lemma}
Soient $k$ un corps de caractéristique $p > 0$, $X \to \Spec k$ un
morphisme de type fini et $G$ un groupe fini agissant sur $X \to \Spec
k$. Supposons que $G$ possède un unique sous-groupe de $p$-Sylow $H$
(en particulier, $H$ est distingué dans $G$) alors on a un
isomorphisme canonique pour tout $i \ge 0$ 
$$\Ext^i_G(\Omega_{X/k}, \O_X) \to \Ext^i_H(\Omega_{X/k}, \O_X)^{G/H}.$$
\end{lemma}

\preuve On a une suite exacte
$$\textrm{III}_2^{j, \ell}=\H^j(G/H, \Ext^\ell_H(\Omega_{X/k},
\O_X))$$ qui converge vers $\Ext^{j+\ell}_G(\Omega_{X/k}, \O_X)$. Or
le cardinal de $G/H$ est inversible dans $k$ donc $\textrm{III}_2^{j,
\ell}=0$ pour tout $j \ge 1$. Le résultat est donc immédiat. \qed 

Soit $C \to \Spec k$ une courbe semi-stable sur un corps
algébriquement clos et munie de l'action d'un groupe fini $G$. Alors
on a, comme dans \cite{BertinMezard}, un morphisme local-global
reliant les déformations équivariantes globales du premier ordre aux
déformations équivariantes <<locales>> du premier ordre et s'insérant
dans une suite exacte 
\begin{multline*}
0 \to \H^1(C/G, \ext^0_G(\Omega_{C/k}, \O_C)) \to
\Ext^1_G(\Omega_{C/k}, \O_C) \to  \\ \to \H^0(C/G,
\ext^1_G(\Omega_{C/k}, \O_C))\to \H^2(C/G, \ext^0_G(\Omega_{C/k},
\O_C)). 
\end{multline*}
D'autre part, comme $C$ est une
courbe, on a $\H^2(C/G, \ext^0_G(\Omega_{C/k}, \O_C))=0$, ce qui donne
une première filtration de $\Ext^1_G(\Omega_{C/k}, \O_C)$. Une
filtration de ce type a déjà été utilisée dans \cite{BertinMezard}
afin d'obtenir la surjectivité du morphisme local-global. Toutefois,
dans le cas des courbes non irréductibles, le faisceau
$\ext^1_G(\Omega_{C/k}, \O_C)$ n'est pas, en général, à support dans
un sous-schéma de dimension zéro (ce sera le cas si l'action est libre
sur un ouvert dense) et ne fournit donc pas de morphisme
local-global. 

\section{Étude des déformations locales}

Nous nous intéressons dans cette section aux 
déformations locales des courbes semi-stables. 

\subsection{Déformation des points lisses}
\label{DeformationPointLisse}
Soit $k$ un corps algébriquement clos de caractéristique $p>0$. Notons
$R=k[[x]]$, $X=\Spec R$, $\Omega=k[[x]] dx$ le module des
différentielles complété et $\nu_x$ la valuation normalisée de $R$.  

L'une des notions fondamentales dans l'étude  des déformations
équivariantes est celle de conducteur, et plus généralement celle de
suite de sauts de ramification.

\begin{definition}
\label{DefConducteurLisse}
Soit $G$ un $p$-groupe fini agissant sur $X=\Spec R$. Notons $G_0$
l'image de $G$ dans $\Aut X$ et  
$$G_i:=\{\sigma \in G_0 | \nu_x(\sigma(x)-x) \ge i+1\}$$ 
 les groupes de ramification supérieurs. La suite croissante $m_0,
\ldots, m_{n-1}$ des entiers $j \in \NN$ tels que $G_j \not = G_{j+1}$
sera appelée la \textit{suite des sauts de ramification}. Si $G_0 \not
= \{Id\}$, l'entier $m_0$ sera appelé le \textit{conducteur} de $G$
(pour son action sur $R$). Dans le cas contraire, le conducteur de $G$
sera par convention $\infty$. 
\end{definition}

On remarquera que, dans la définition précédente, on ne suppose pas
que le groupe $G$ agit librement sur le corps des fractions. Les
groupes de ramification considérés sont donc ceux de l'image de $G$
dans $\Aut X$. 

\begin{proposition}
\label{ProprieteConducteur}
Soit $G$ un $p$-groupe fini agissant sur $X=\Spec R$. Notons $G_0$
l'image de $G$ dans $\Aut X$ et supposons le non trivial. Alors la
suite $m_0, \ldots, m_{n-1}$ des sauts de ramification vérifie les
propriétés suivantes : 
\begin{enumerate}[{\rm(i)}]
\item Pour tout $j \in \{1, \ldots, n-1\}$ on a $m_j = m_0 \mod p$.
\item On a $p \nmid m_0$.
\item Si $G$ est cyclique alors $|G_0| = p^{n}$.
\end{enumerate}
\end{proposition}

\preuve (i) cf. \cite{Corps_Locaux} Proposition IV.2.11. 
(ii) Soit $\sigma\in G$ d'ordre $p$. La théorie d'Artin-Schreier
permet de voir que $m:=\nu_x(\sigma(x)-x)-1$ est premier à $p$. Le
résultat découle alors du point (i) car $m$ est un saut de la
filtration de ramification de $G$. 
(iii) Les groupes $G_i$ sont cycliques car $G$ est cyclique et comme
$G_i/G_{i+1}$ est un produit direct de groupes cycliques d'ordre $p$
d'après \cite{Corps_Locaux} Corollaire IV.2.3, on a $G_i/G_{i+1}=0$ ou
$G_i/G_{i+1}=\ZZ/p\ZZ$. Le résultat est donc immédiat. \qed

\begin{exemple}
\label{ActionpCyclique}
Soient $k$ un corps de caractéristique $p$ et $m$ un entier premier à $p$. Soit $\sigma$ un générateur de $\ZZ/p\ZZ$. Définissons une action à gauche de $\ZZ/p\ZZ$ sur $k[[x]]$ par
$$\sigma^i(x)= {x \over (1+i x^m)^{1/m}}.$$
On montre aisément que ceci définit une action de $\ZZ/p\ZZ$ sur $k[[x]]$ de conducteur $m$. 
Notons $\d x$ la base de $\Hom_R(\Omega, R)$ telle que $\d x (dx)=1$. On a alors 
$$\sigma^i(dx)={1 \over (1+ix^m)^{1+1/m}} dx, \ \sigma^i\left(\d x\right) = (1+i x^m)^{1+1/m}\d x$$ 
(car $\left(\sigma^i\left(\d x\right)\right)(dx)=\sigma^i\left(\d x (\sigma^{-i}(dx))\right )$).
D'autre part on peut montrer que, quitte à changer d'uniformisante
$x$, toute action fidèle de $\ZZ/p\ZZ$ sur $k[[x]]$ ayant un
conducteur $m$ est de cette forme (cf. \cite{BertinMezard}, Lemme
4.2.1).  
\end{exemple}

Pour tout nombre réel $\alpha$, notons $\lfloor \alpha \rfloor$ la
partie entière inférieure de $\alpha$ et $\lceil \alpha \rceil$ la
partie entière supérieure. On peut alors calculer la dimension de
l'espace tangent au foncteur des déformations équivariantes dans le
cas cyclique. 

\begin{theorem}[\cite{BertinMezard}, Théorème 4.1.1]
\label{DimensionLisse}
Soient $G$ un groupe cyclique d'ordre $p^n$ agissant fidèlement sur
$\Spec R$ et $\dif$ la valuation de la différente de l'extension
$R/R^G$.  
On a 
$$\dim_k \Ext^1_G(\Omega, R)=\left\lfloor { 2 \dif \over p^n } \right\rfloor-\left\lceil { \dif \over p^n } \right\rceil.$$
\end{theorem}

Nous aurons besoin par la suite de savoir s'il existe des éléments de
$\Ext^1_G(\Omega, R)$ possédant certaines propriétés. Plus
précisément, nous aurons besoin du résultat suivant.

\begin{theorem}
\label{ElementNonNulTraceNulle}
Soient $G$ un $p$-groupe cyclique agissant sur $X=\Spec R$, $G_0$
l'image de $G$ dans $\Aut X$ et $\dif$ la valuation de la différente de l'extension
$R/R^G$. Alors il existe une base $\phi$ du $R$-module $\Hom(\Omega,
R)$ telle que $\Tr_G \phi = 0$ si et seulement si $2\dif+1 \not = 0
\mod |G_0|$ ou si $G \not = G_0$. 
\end{theorem}

Avant de passer à la démonstration du théorème, nous
avons besoin de quelques résultats préliminaires.

\begin{lemma}
\label{TraceAdditive}
Soient $G \subset \Aut_k(R)$ un groupe fini et $H \subset G$ un
sous-groupe distingué. Notons $y:= \prod_{\sigma \in H}
\sigma(x)$. Pour tout $\phi \in \Hom_{R}(\Omega,R)$ on a $$\Tr_G
\phi:=\sum_{\sigma \in G} \sigma. \phi = \Tr_{G/H} \left ( \left (
\Tr_{H} \phi \right)(dy)\d y \right ).$$  
\end{lemma}

\preuve C'est un corollaire immédiat de la multiplicativité de la
différente (qui permet de voir que $\d x=\d x(dy) \d y$)  et du fait
que $\Tr_G = \Tr_{G/H} \circ \Tr_H$. \qed

\begin{lemma}
\label{CalculValuation}
Soient $G$ un groupe cyclique d'ordre $p^n$ agissant fidèlement sur
$R$ et $m_0, \ldots, m_{n-1}$ la suite des sauts de ramification. 
Soient $f \in R$, $\ell:=\nu_x(f)$ et  $z:=\prod_{\tau \in G}
\tau(x)$. Supposons que $2m_0+1 = 0 \mod p$ et que  
$$\ell+\sum_{i=0}^{n-2} p^i(2m_{n-i-1}+1)=0 \mod p^n$$
alors 
$$\nu_z \left ( \left ( \Tr_G f \d x \right ) (dz) \right )={1 \over
p^n}\left ( \ell+(p-1) \left ( \sum_{i=0}^{n-1} (2m_{n-1-i}+1)p^i
\right )\right ).$$  
\end{lemma}

\preuve On voit que cette propriété ne dépend pas de l'uniformisante
$x$, on pourra donc s'autoriser à en changer.

Nous allons montrer ce lemme par récurrence sur $n$.
Supposons $n=1$, notons $\sigma$ un générateur de $G$ et écrivons
$f=\sum_{i=\ell}^\infty a_i x^i$. Les hypothèses du lemme impliquent que
$p | \ell$. 

Quitte à changer d'uniformisante, on peut supposer (cf. l'exemple
\ref{ActionpCyclique}) que l'action de $G$ est donnée par  
$$\sigma(x)={x \over (1+x^{m_0})^{1/m_0}}.$$
On a donc $$\left ( \Tr_G  f \d x \right ) (dx) = \sum_{q=0}^{p-1}
\sum_{i=\ell}^\infty a_i x^i (1+qx^{m_0})^{1+1/m_0-i/m_0}.$$ 
Notons $\binome{\alpha}{n}$ les coefficients binômiaux.
On a
\begin{eqnarray*}
\sum_{q=0}^{p-1}(1+qx^{m_0})^{1+1/m_0-i/m_0} & = & \sum_{j=0}^\infty
\binome{1+1/m_0-i/m_0}{j}x^{m_0j}
\underbrace{\sum_{q=0}^{p-1}q^{j}}_{\begin{tiny}
=\begin{cases}0 \textrm{ si } p-1 | j \\
1 \textrm{ sinon}\end{cases}\end{tiny}} \\ 
& = & -\sum_{j'=1}^\infty \binome{1+1/m_0-i/m_0}{(p-1)j'}x^{m_0(p-1)j'}.
\end{eqnarray*}

Le terme de degré $m_0(p-1)+\ell$ de $\left ( \Tr_G f \d x \right )(dx)$ est donc
$$a_\ell \binome{1+1/m_0-\ell/m_0}{(p-1)}.$$ Ce terme est non nul car
$a_{\ell}\not=0$ et le terme binômial est non nul car $\ell=2m_{0}+1=0\mod p$.

En écrivant $dz=\d{x} (dz)dx$ et en utilisant la $k[[x]]$-linéarité de
$\Tr_{G} (f \d{x})$ on déduit
\begin{eqnarray*}
\nu_x \left ( \left ( \Tr_G f \d x\right )(dz) \right )& = & \nu_x
\left ( \d x(dz)\right ) + \nu_x\left ( \left (\Tr_G  f\d x\right )
(dx)\right) \\ 
& = & (m_0+1)(p-1)+\ell+m_0(p-1).
\end{eqnarray*}
On a donc $\nu_z\left ( \left ( \Tr_G f \d x\right )(dz) \right ) = {1
\over p}(\ell+(2m_0+1)(p-1))$ car $k[[x]]/k[[z]]$ est totalement ramifiée
de degré $p$. 

\medskip

Supposons le résultat démontré au rang $n$ et montrons-le au rang
$n+1$. On a alors $$\ell+(p-1) \sum_{i=0}^{n-1} p^i (2m_{n-i}+1)=0 \mod
p^{n+1}$$ 
et $2m_n+1=0 \mod p$ car $m_n=m_0 \mod p$ (cf. proposition \ref{ProprieteConducteur}).
Notons $H$ un sous-groupe de $G$ d'ordre $p$ et $y:=\prod_{\sigma \in
H} \sigma(x)$. D'après le calcul effectué pour $n=1$ on trouve  
$$\nu_y\left ( \left (\Tr_H f \d x\right) (dy) \right)={1 \over
p}(\ell+(2m_n+1)(p-1))$$ (c'est bien un entier car $\ell+(p-1)(2m_n+1)=0
\mod p$). 

Notons $Q:=\left (\Tr_H f {\partial \over \partial x} \right ) (dy)$
et $\ell':=\nu_y(Q)$, on a  
$$\ell'+(p-1) \sum_{i=0}^{n-2} p^{i} (2m_{n-1-i}+1)=0 \mod p^{n}.$$    

On peut donc appliquer l'hypothèse de récurrence à $Q \d y$ et on trouve 
$$\nu_z \left ( \left ( \Tr_{G/H} Q \d y \right )(dz) \right )={1
\over p^n}\left (\nu_y(Q)+(p-1) \sum_{i=0}^{n-1} p^i (2m_{n-i-1}+1)
\right )$$ car les sauts de ramification pour l'action de $G/H$ sur
$k[[y]]$ sont donnés par la suite  $\{m_0, \ldots, m_{n-1}\}$ d'après
\cite{Corps_Locaux} Corollaire IV.1. 
Le lemme \ref{TraceAdditive} permet d'écrire 
\begin{multline*}
\nu_z\left ( \left ( \Tr_G f \d x \right ) (dz) \right )  = \nu_z\left
( \left ( \Tr_{G/H} Q\d y \right )(dz) \right ) \\ 
 =  {1 \over p^n}\left ({1 \over p}(\ell+(p-1)(2m_n+1))+(p-1)
\sum_{i=0}^{n-1} p^i (2m_{n-i-1}+1) \right ) 
\end{multline*}
ce qui démontre le résultat. \qed

\begin{lemma}
\label{ElementTraceNulleValuationFixee}
Soit $\sigma$ un automorphisme d'ordre $p$ de $k[[x]]$ de conducteur
$m<\infty$. Supposons que $2m+1 = 0 \mod p$ et donnons-nous $q \in \NN$
premier à $p$. Alors il existe $f \in k[[x]]$ tel que $\nu_x(f)=q$ et
$\Tr_{\langle \sigma \rangle} f {\partial \over \partial x} = 0$.  
\end{lemma}

\preuve Notons $z:=\prod_{i=0}^{p-1} \sigma^i(x)$ et 
$$\theta: \Hom_{k[[x]]}(k[[x]]dx, k[[x]]) \rightarrow
k[[z]]=k[[x]]^{\langle \sigma \rangle}$$ 
l'application définie par $\theta(\phi):=\left ( \Tr_{\langle \sigma
\rangle} \phi \right)(dz)$ 
(elle est bien à valeurs dans $k[[z]]$ car $\theta(\phi)$ est fixe sous $\sigma$).
L'application $\theta$ est un morphisme de $k[[z]]$-modules et son
image est un sous-$k[[z]]$-module de $k[[z]]$, donc $\Im \theta$ est
de la forme $z^{r} k[[z]]$. Montrons, par l'absurde, que $r={1 \over p}(2m+1)(p-1)$.  

D'après le lemme \ref{CalculValuation} (cas $n=1$) on a ${1 \over
p}(2m+1)(p-1)=\nu_z\left ( \theta(\d x) \right ) $ et donc $r < {1
\over p}(2m+1)(p-1)$. Soit $h \in k[[x]]$ telle que  $\nu_z(\theta(h
\d x))=r$. On a nécessairement $\nu_x(h) > 0$ car sinon le lemme
\ref{CalculValuation} imposerait $r=\nu_z(\theta(h \d x))={1 \over
p}(2m+1)(p-1)$ ce qui est absurde. Or $h+1$ est non nulle en zéro donc
on a $\nu_z(\theta((h+1)\d x)))={1 \over p}(2m+1)(p-1)$ d'après le
lemme \ref{CalculValuation}. Mais on a  
$$\nu_z\left(\theta\left((h+1)\d
x\right)\right)=\nu_z\left(\theta\left(h \d x\right)+\theta\left(\d
x\right)\right)=\nu_z\left(\theta\left(h \d x\right)\right)=r$$ 
car $\nu_z(\theta(h \d x)) = r < {1 \over p}(2m+1)(p-1)=\nu_z(\theta(
\d x))$. On a donc bien $r={1 \over p}(2m+1)(p-1)$. 

\medskip

Par suite, notons $q=q'+\ell p$ avec $1 \le q' \le p-1$. Posons
$h:=\theta(x^{q'} \d x)\in k[[z]]$. Alors $g:=h/\theta(\d x )\in
k[[z]]$ car $\theta\left( \d x\right)$ engendre l'image de $\theta$
(en tant que $k[[z]]$-module). Notons $f:=z^\ell(x^{q'}-g)\in k[[x]]$. On
a alors  
$$\theta\left (f\d x\right )=\theta\left (z^\ell (x^{q'}-g)\d x\right
)=z^\ell \left( h-g\theta\left(\d x\right)\right)=0,$$ 
c'est-à-dire $\Tr_{\langle \sigma \rangle} f \d x=0$ (car l'évaluation en $dz$ est injective). 
De plus, cela prouve que $\nu_z(g) > 0$ car dans le cas contraire on
aurait $\nu_x(x^{q'}-g)=0$ ce qui impose, par le lemme
\ref{CalculValuation}, que $$\nu_z\left(\theta\left((x^{q'}-g)\d
x\right)\right) = {1 \over p} (p-1)(2m+1) < \infty$$ ce qui est
absurde. On a donc $\nu_x(g) \ge p$ et  
$$\nu_x(f)=\nu_x(z^\ell (x^{q'}-g))=\nu_x(z^\ell)+\nu_x(x^{q'}-g)=p\ell+\nu_x(x^{q'})=q$$ car $q' < p \le \nu_x(g)$. \qed

\noindent \textit{Preuve du théorème \ref{ElementNonNulTraceNulle} :} Si $G_0 \not = G$ alors le résultat est trivial car $\Tr_G=0$. On peut donc supposer que $G=G_0$ (c'est-à-dire que l'action est fidèle).

Notons $m_0, \ldots, m_{n-1}$ la suite des sauts de ramification. On a alors 
$$\dif = (p-1)\sum_{i=0}^{n-1} p^i(m_{n-i-1}+1)$$
 (car $\dif=\sum_i \card(G_i)-1$, cf. par exemple \cite{Corps_Locaux}
Chap IV, \S 2, Proposition 4) et

$$2 \dif + 1=p^n+(p-1)\sum_{i=0}^{n-1} p^i(2m_{n-i-1}+1).$$ 
C'est sous cette forme que nous allons utiliser la différente. 

Si $2\dif+1 = 0 \mod p^n$, donc $\sum_{i=0}^{n-1} p^i(2m_{n-i-1}+1)=0
\mod p^n$ (en particulier on a $2m_0+1=0 \mod p$), alors d'après le
lemme \ref{CalculValuation}, pour tout $f \in k[[x]]$ inversible on a

$$\nu_z\left ( \left ( \Tr_G f {\partial  \over \partial x} \right )
(dz) \right )={1 \over p^n}\left ( (p-1) \left ( \sum_{i=0}^{n-1}
(2m_{n-1-i}+1)p^i \right )\right )< \infty,$$

et donc $\Tr_G\left( f \d x \right ) \not = 0$.
\medskip

Supposons maintenant que $\sum_{i=0}^{n-1} p^i(2m_{n-i-1}+1)\not =0
\mod p^n$.  Nous allons exhiber $\phi$ comme dans l'énoncé du
théorème. Dans un premier temps, nous allons supposer de plus que
$2m_0+1 \not = 0 \mod p$. Soit $\sigma \in G$ un élément d'ordre
$p$. Si $m_{n-1} = -1 \mod p$ alors on pose $b := -1$ et $a :=
(1+m_{n-1})/p \in \NN$ de sorte que $ap+bm_{n-1}=1$. Si $m_{n-1} \not
= -1 \mod p$ alors il existe $a, b \in \ZZ$ tel que $ap+bm_{n-1}=1$ et
$1 \le b < p-2$ (car $2m_{n-1}+1 \not = 0 \mod p$). Quitte à changer
d'uniformisante $x$, on peut supposer que l'action de $\sigma$ est
donnée par

$$\sigma(x)={x \over (1+x^{m_{n-1}})^{1/m_{n-1}}}.$$ 
 Posons 
$$\phi:=\left(-1+x^{m_{n-1}(p-1)}\right)^{a/m_{n-1}} \d x.$$ 
On a
\begin{eqnarray*}
\left(\Tr_{\langle \sigma \rangle} \phi\right)(dx) & = &
\sum_{q=0}^{p-1} \left (-1+{x^{m_{n-1}(p-1)} \over (1+qx^{m_{n-1}})^{p-1}}\right )^{a/m_{n-1}} (1+qx^{m_{n-1}})^{1/m_{n-1}+1} \\
& = & \left (-1+x^{m_{n-1}(p-1)}\right )^{a/m_{n-1}} \sum_{q=0}^{p-1} (1+qx^{m_{n-1}})^{(b+1)} \\
& = & \left (-1+x^{m_{n-1}(p-1)}\right )^{a/m_{n-1}} \sum_{i=0}^{b+1}
\binome{b+1}{i} x^{i m_{n-1}} \underbrace{\sum_{q=0}^{p-1} q^i.}_{=0
\textrm{ car } i\le b+1 < p-1}
\end{eqnarray*}
Or 
$$\Tr_G \phi = \Tr_{G/\langle \sigma \rangle} \left ( \Tr_{\langle
\sigma \rangle} \phi \right ) = 0$$ et comme $\phi$ est une base de
$\Hom(\Omega, R)$, on a prouvé le lemme dans le cas $2m_0+1=0 \mod
p$. 

\medskip
Supposons finalement que $\sum_{i=0}^{n-1} p^i(2m_{n-i-1}+1)\not =0
\mod p^n$  et que $2m_0+1 = 0 \mod p$. En particulier on a
$p^{n-1}(2m_0+1)= 0 \mod p^n$. Notons  
$$n_0:=\min\left\{1 \le j \le n-1 | \sum_{i=0}^{j-1} p^i(2m_{n-i-1}+1)
\not = 0 \mod p^{j+1}\right\},$$ 
l'ensemble de droite est non vide car il contient $n-1$. On a $n_0 >
1$ car $2m_{n-1}+1 =2m_0+1  = 0 \mod p$. Par minimalité de $n_0$ (et
comme $n_0 > 1$) on a  
$$\sum_{i=0}^{n_0-2} p^i(2m_{n-i-1}+1) = 0 \mod p^{n_0}.$$
Notons $y:=\prod_{\sigma \in G_{n_0}} \sigma(x)$. Le lemme
\ref{CalculValuation} nous donne  
$$q:=\nu_y\left ( \left ( \Tr_{G_{n_0}} \d x \right ) (dy) \right )={1
\over p^{n_0}} \left ( (p-1)\sum_{i=0}^{n_0-1} p^i(2m_{n_0-1-i}+1)
\right )$$ 
et donc  $q\not = 0 \mod p$ par définition de $n_{0}$.
Soit $\sigma$ un automorphisme d'ordre $p$ de $G/G_{n_0}$ (il en
existe car $n_0 > 1$). Notons $h:=\left ( \Tr_{G_{n_0}} \d x \right )
(dy)$ et $f\in k[[y]]$ tel que $\nu_y(f)=q$ et $\Tr_{\langle \sigma
\rangle} f \d y = 0$ (il en existe d'après le lemme
\ref{ElementTraceNulleValuationFixee}). Comme $h$ et $f$ ont même
valuation, $f/h$ est un inversible de $k[[y]]$ donc, en appliquant le
lemme \ref{TraceAdditive}, on a 
\begin{eqnarray*}
\Tr_{G_{n_0-1}} f/h \d x & = & \Tr_{\langle \sigma \rangle} \left (
\left( \Tr_{G_{n_0}} f/h \d x\right)(dy) \d y \right ) \\ 
& \hskip -2cm= & \hskip -1cm\Tr_{\langle \sigma \rangle} \left ( f/h
\left( \Tr_{G_{n_0}} \d x \right)(dy) \d y \right ) \textrm{ car } f/h
\in k[[x]]^{G_{n_0}}=k[[y]] \\ 
& = & \Tr_{\langle \sigma \rangle} \left ( f {\partial \over \partial y} \right ) \textrm{ par définition de  } h \\
& = & 0 \textrm{ par définition de } f.
\end{eqnarray*}
Par suite, on trouve que $\Tr_G\left ( f/h \d x \right )=0$. \qed

Si on suppose $G$ cyclique, les éléments $\phi \in \Hom_R(\Omega, R)$
tels que $\Tr_G \phi = 0$ possèdent une interprétation simple.  
En effet, si $G$ est engendré par $\sigma$ alors on a l'égalité
$$\H^1(G, \Hom_R(\Omega, R))=\{ \phi \in \Hom_R(\Omega, R) | \Tr_G
\phi = 0 \} / \{\phi-\sigma \phi\}_\phi.$$ Notons $\X=\Spec
k[\epsilon]/(\epsilon^2)[[x]]$ et choisissons $\phi \in \Hom(\Omega,
R)$ tel que $\Tr_G \phi=0$. La déformation équivariante définie par
$\phi$ est donnée par 
$$\sigma_{\phi}^i(x):=\sigma^i(x)+\epsilon \sum_{j=1}^i \sigma^j \left(\phi \right )(dx).$$
Par suite, si $\phi(dx)\not\in xk[[x]]$ alors $\sigma_\phi(x)\not\in
xk[\epsilon]/ (\epsilon^2) [[x]]$. D'autre part, si on a une
déformation équivariante définie par $\phi\in \H^1(G, \Hom_R(\Omega,
R))$ telle que $\sigma_\phi(x) \not\in xk[\epsilon]/ (\epsilon^2)
[[x]]$ alors $\phi(dx) \not \in xk[[x]]$. 

L'existence d'une base $\phi\in\Hom(\Omega, R)$ telle que $\Tr_G
\phi=0$ sera particulièrement importante quand il s'agira d'étudier
les déformations équivariantes du point double. 

\subsection{Déformation des points doubles}

Soit $k$ un corps algébriquement clos de caractéristique $p$.
Considérons l'anneau $R=k[[x, y]]/(xy)$ et notons $R'=k[[x]]\times
k[[y]]$ sa normalisée, $X=\Spec R$, $X'=\Spec R'$, $\Omega$ le
complété de $\Omega_{X/k}$ et $\Omega'$ le complété de
$\Omega_{X'/k}$. 
 
\begin{definition}
Reprenons les notations ci-dessus.
Les composantes connexes de $X'$ seront appelées \textit{branches} de $X$. 
\end{definition}

\begin{definition}
Soit $G$ un $p$-groupe agissant sur $X$. Notons $H$ le sous-groupe de
$G$ formé des éléments ne permutant pas les branches. Notons $m$ et
$m'$ les conducteurs de $H$ sur chaque branche de $X$ (voir la
définition \ref{DefConducteurLisse}), $\dif$ et $\dif'$ les
différentes. Le couple $(m, m')$ sera appelé \textit{conducteur} de
$G$ sur $X$ et le couple $(\dif, \dif')$ la \textit{différente} de $G$
sur $X$. 
\end{definition}

On peut remarquer que la notation précédente fait référence à un ordre
sur les branches. Dans les faits, la première coordonnée de chaque
couple sera toujours relative à l'action sur $x$ et la deuxième à
l'action sur $y$. 
 
\begin{exemple}
\label{ActionpCycliqueSingulier}
Soient $m, m'\in \NN$ deux entiers premiers à $p$ et $\sigma$ un
générateur de $\ZZ/p\ZZ$. Définissons une action de $\ZZ/p\ZZ$ sur
$k[[x, y]]/(xy)$ par  
$\sigma^i(x)={x \over (1+ix^m)^{1/m}}$ et $\sigma^i(y)={y \over
(1+iy^{m'})^{1/m'}}$. On montre aisément que ceci définit une action
de $\ZZ/p\ZZ$ sur $X$ de conducteur $(m, m')$. De la même manière que
dans l'exemple \ref{ActionpCyclique}, on peut montrer que, quitte à
changer d'uniformisante, toutes les actions de $\ZZ/p\ZZ$ sur $X$ qui
sont fidèles sur chacune des branches de $X$ sont de cette forme.  
\end{exemple}

Nous aurons besoin dans la suite des expressions des automorphismes en
fonction des paramètres uniformisants, ceci est donné par le lemme
suivant qui fixera également les notations. On prendra garde au fait
que l'action à gauche sur le schéma donne une action à droite sur
l'anneau de fonctions. Toutefois l'automorphisme $\sigma$ du schéma
sera encore noté $\sigma$ lorsqu'il agira sur l'anneau de fonctions
(pour ne pas alourdir les notations). 

\begin{lemma}
\label{CoordonneeAutomorphisme}
Soit $\sigma$ un $k$-automorphisme de $R$ qui fixe les branches. Alors
il existe des séries formelles uniques $P_0 \in xk[[x]]$ et $P_1 \in
yk[[y]]$ telles que $\sigma$ soit défini par  
$$\sigma(x)=P_0(x)\hskip 1cm \sigma(y)=P_1(y).$$
\end{lemma}

\preuve Comme $\sigma$ n'échange pas les branches, l'image de $x$ est
dans l'idéal premier engendré par $x$ et comme $xy=0$, on obtient le
résultat. L'unicité est immédiate. \qed 

Pour un automorphisme $\sigma$ de $X$ ne permutant pas les branches,
le lemme précédent permet de considérer $\sigma(x)$
(resp. $\sigma(y)$) comme une série formelle en $x$ (resp. en $y$). 
Si $\sigma$ est d'ordre $n$ on voit que le terme de degré $1$
de la série formelle $\sigma (x)$ est une racine $n$-ième de l'unité.

Nous allons répartir les déformations équivariantes en deux classes de
nature différente. Pour cela, souvenons-nous que toute déformation
(non équivariante) de $X$ au-dessus d'un anneau artinien $A$ dont le
corps résiduel est $k$ est de la forme $\Spec A[[x, y]]/(xy-a)$ avec
$a \in A$ (cf. \cite{Deligne_Mumford}). Il est donc naturel de poser
la définition suivante. 

\begin{definition}
\label{DefTopTriv}
Soient $A$ un objet de $\A$ et $(\X\to\Spec A, G)$ une déformation
équivariante de $(X, G)$. On dira que $(\X, G)$ est une déformation
\textit{topologiquement triviale} si $\X \cong \Spec A[[x, y]]/(xy)$
en tant que schéma.  
\end{definition}

La suite spectrale 
$\textrm{II}_2^{j,\ell}:=\H^j(G, \Ext^{\ell} (\Omega, R))$ fournit une
suite exacte
\begin{equation}
\label{EquationPointDouble}
0 \rightarrow \H^1(G, \Ext^0(\Omega, R)) \rightarrow \Ext^1_G(\Omega,
R) \overset{{\textit{ép}}}{\rightarrow} \H^0(G, \Ext^1(\Omega, R)). 
\end{equation}
Les déformations équivariantes du premier ordre dont l'image dans
$\Ext^1(\Omega, R)$ est nulle sont les déformations topologiquement
triviales, celles-ci s'identifient donc avec le groupe $\H^1(G,
\Ext^0(\Omega, R))$.  

Nous allons maintenant décrire l'espace $\H^1(G, \Ext^0(\Omega, R))$
en terme des déformations de la normalisée. Ensuite, nous calculerons
la dimension de $\Im {\textit{ép}}$. 

\subsubsection{Lien avec la normalisée}
\label{LienNormaliseeLocal}

Soient $X=\Spec k[[x, y]]/(xy)$ et $G$ un $p$-groupe agissant sur $X$. 
Nous allons décomposer l'espace $\H^1(G, \Ext^0(\Omega, R))$ grâce à
la suite exacte \eqref{SuiteExacteDeformationToNormalisee} et aux
propositions \ref{LocalSingulier1} et \ref{LocalSingulier2}.  

\begin{lemma}
On a une suite exacte canonique
\begin{equation}
\label{MorphismeNormaliseeLocal}
0 \to \Hom_R(\Omega, R) \to \Hom_{R'}(\Omega', R') \to \Hom_R(\Omega, R'/R) \to 0.
\end{equation}
Notons $\d{x}$ (resp. $\d{y}$) l'élément de $\Hom_{R'}(\Omega',
R')$ tel que $\d{x} (dx)=1$ et $\d{x} (dy)=0$ (resp. $\d{y} (dy)=1$ et
$\d{y} (dx)=0$). Alors on a une identification naturelle 
\begin{equation}\label{IdentificationDualOmega}
\Hom_R(\Omega, R)\cong xk[[x]]\d{x}\oplus yk[[y]]\d{y}.
\end{equation}
\end{lemma}

\preuve On montre aisément qu'on a un isomorphisme canonique
$$\Hom_R(\Omega', R') \rightarrow \Hom_{R}(\Omega, R').$$
D'autre part, on a une suite exacte de $R[G]$-modules
$$0 \to R \to R' \to R'/R\to 0.$$
En appliquant le foncteur $\Hom_R(\Omega, \bullet)$ à cette suite
exacte, on obtient une suite exacte longue
$$0 \to \Hom_R(\Omega, R) \to \Hom_{R'}(\Omega', R') \to
\Hom_R(\Omega, R'/R) \to \ldots$$
Il suffit donc de montrer que la dernière flèche est surjective. Or 
on voit que l'espace vectoriel $\Hom_{R}(\Omega, R'/R)$ est engendré
par les images de $\d x$ et $\d y$. \qed 

La description précédente permet, dans une large mesure, de se ramener
au cas où $G$ agit sur $X$ sans permuter les branches. En effet,
on a le corollaire suivant.

\begin{corollary}
\label{PermutteToNonPermutte}
Soient $G$ un groupe fini agissant sur $X$ et $H$ le sous-groupe de
$G$ constitué des éléments ne permutant pas les branches de $X$. Alors
$H$ est distingué dans $G$ et le morphisme canonique 
$$\Ext^1_G(\Omega, R) \to \Ext^1_H(\Omega, R)^{G/H}$$ est un
isomorphisme. 
\end{corollary}

\preuve Le fait que $H$ est un sous-groupe distingué est évident. La
suite spectrale $\textrm{III}_{2}^{j, \ell}:=\H^j(G/H, \Ext^{\ell} (\Omega,
R))$ fournit une suite exacte 
\begin{multline*}
0 \to \H^1(G/H, \Ext^0_H(\Omega, R)) \to \Ext^1_G(\Omega, R) \to
\H^0(G/H, \Ext^1_H(\Omega, R)) \to  \\ \H^2(G/H, \Ext^0_H(\Omega,
R)). 
\end{multline*}
Il est alors aisé de voir, grâce à l'isomorphisme
\eqref{IdentificationDualOmega}, que $\Ext^0_H(\Omega, R)$ 
est un $k[G/H]$-module libre car $G/H$ échange les branches
analytiques de $X/H$ (et ne fait que ça) donc pour $j \ge 1$ on a 
$$\H^j(G/H, \Ext^0_H(\Omega, R))=0.$$  \qed 

Nous nous limiterons donc désormais au cas où les groupes agissent
sans permuter les branches (ce qui est le cas général quand on se
restreint aux actions des $p$-groupes).

En considérant les éléments fixes sous l'action de $G$ dans la suite
exacte \eqref{MorphismeNormaliseeLocal}, on obtient une
suite exacte
\begin{multline*}
0 \to \Hom_R(\Omega, R)^G \to \Hom_{R'}(\Omega', R')^G \to
\Hom_R(\Omega, R'/R)^G  \to \\ 
\to \H^1(G, \Hom_R(\Omega, R)) \to \H^1(G, \Hom_{R'}(\Omega', R')) \to
\H^1(G, \Hom_R(\Omega, R'/R)). 
\end{multline*}

Notons $\phi_x$ et $\phi_y$ les $1$-cocycles $G \to \Hom_R(\Omega, R)$
définis par $\phi_x(\sigma)=\d x-\sigma \d x$ et $\phi_y(\sigma)=\d
y-\sigma \d y$ (dont on montre qu'ils sont effectivement à valeurs
dans $\Hom_R(\Omega, R)$).  
On voit alors que le morphisme
$$\Hom_R(\Omega, R'/R)^G  \to \H^1(G, \Hom_R(\Omega, R))$$
se factorise par  
$$\left ( k \phi_x \oplus k \phi_y \right )^G\to \H^1(G, \Hom_R(\Omega, R)).$$
Finalement, on a une suite exacte
\begin{multline}
\label{SuiteExacteDeformationToNormalisee}
\left(k \phi_x \oplus k \phi_y\right)^G \to \H^1(G, \Hom_R(\Omega, R))
\to \H^1(G, \Hom_{R'}(\Omega', R')) \to \\ \H^1(G, \Hom_R(\Omega,
R'/R)). 
\end{multline}
Nous allons maintenant préciser le noyau du morphisme
$$\H^1(G, \Hom_R(\Omega, R)) \to \H^1(G, \Hom_{R'}(\Omega', R')).$$

\begin{proposition}
\label{LocalSingulier1}
Soit $G$ un $p$-groupe agissant sur $X$ avec conducteur $(m, m')$ et
sans permuter les branches. 
\begin{enumerate}[{\rm(i)}]
\item Si $m < \infty$ et $m+1 \not = 0 \mod p$ (resp. $m'< \infty$ et
$m'+1 \not = 0 \mod p$), alors $\phi_x$ (resp. $\phi_y$) est non nul
dans $\H^1(G,\Hom_R(\Omega, R))$.  

De plus, si $m < \infty$, $m+1 \not = 0 \mod p$, $m'< \infty$  et
$m'+1 \not = 0 \mod p$ alors les images de $\phi_x$ et de $\phi_y$
dans $\H^1(G, \Hom_R(\Omega, R))$ sont linéairement indépendantes. 
\item Si $G$ est un groupe $p$-cyclique et si $m+1=0 \mod p$ ou $m =
\infty$ alors $\phi_x$ est nul dans $\H^1(G, \Hom_R(\Omega, R))$. 
\end{enumerate}
\end{proposition}

\preuve Démontrons tout d'abord le point (i). Supposons que $m <
\infty$ et $m+1 \not = 0 \mod p$. Soit $\sigma\in G$ un élément tel
que $\nu_x(\sigma(x)-x) = m+1$ (sur la normalisée de $X$). Notons
$\sigma(x)=x+a_{m+1}x^{m+1}+\ldots$. On a alors  
$$\phi_x(\sigma)(dx)=1-{1 \over \der x
(\sigma(x))}=1-(1-(m+1)a_{m+1}x^m+\ldots)$$ 
et donc $\nu_x(\phi_x(\sigma)(dx))=\nu_x((m+1)a_{m+1}x^m)=m$.
Si $\phi_x$ est nul dans $$\H^1(G, \Hom_R(\Omega, R))$$ alors il
existe $f \in xk[[x]]$ tel que pour tout $\tau \in G$ on a
$\phi_x(\tau)=f \d x - \tau\left(f \d x \right )$ (car $\Hom_R(\Omega,
R)\cong xk[[x]] \d x \oplus y k[[y]] \d y$, cf. l'isomorphisme
\eqref{IdentificationDualOmega}). 
Soit $n > 0$, on a 
\begin{eqnarray*}
\left( x^n \d x - \sigma\left(x^n \d x \right )\right )(dx) & = & x^n
- {\sigma(x)^n \over \der x \sigma(x)} = {x^n \der x \sigma(x)-\sigma(x)^n \over \der x \sigma(x)} \\
& \hskip -2cm = & \hskip
-1cm{x^n(1+(m+1)a_{m+1}x^m-(1+a_{m+1}x^m+\ldots)^n) \over \der x
\sigma(x)} 
\end{eqnarray*}
Par suite, on trouve que 
$$\nu_x \left( x^n \d x - \sigma\left(x^n \d x \right )\right )(dx) \ge n+m$$
et par linéarité
$$\nu_x\left( \left ( f \d x - \sigma\left(f \d x \right ) \right)(dx)\right) > m$$
ce qui est absurde donc $\phi_x \not = 0$ dans $\H^1(G, \Hom_R(\Omega, R))$.

D'autre part, comme $\Hom_R(\Omega, R) \cong xk[[x]] \d x\oplus y
k[[y]]\d y$ et que $G$ n'échange pas les branches analytiques, si
$\phi_x$ et $\phi_y$ sont non nuls dans $\H^1(G, \Hom(\Omega, R))$
alors ils sont libres. 

Passons à la démonstration du point (ii). Le cas $m=\infty$ est
trivial car alors $\phi_x=0$. Supposons $m+1=0 \mod p$. Quitte à
changer d'uniformisante, on peut supposer que $\sigma(x)={x \over
(1+x^m)^{1/m}}$ (cf. exemple \ref{ActionpCycliqueSingulier}). 
Posons $f:=(1-x^{m(p-1)})^{m+1 \over pm}-1$ (qui est bien définie car $p | m+1$). 
On a 
\begin{eqnarray*}
\sigma^i\left (f \d x\right ) & = & \left (1-{x^{m(p-1)} \over
(1+ix^m)^{p-1}} \right )^{m+1 \over pm} (1+ix^m)^{1+1/m}\d
x-\sigma^i\left( \d x\right ) \\ 
& = & (1-x^{m(p-1)})^{m+1 \over pm}(1+ix^m)^{1-p{m+1 \over pm}+1/m}
\d x -\sigma^i\left( \d x\right ) \\ 
& = &  f \d x + \d x - \sigma^i\left( \d x\right ).
\end{eqnarray*}
On a donc $f\in xR$ et 
$$f \d x - \sigma^i \left (f \d x\right ) = -\phi_x$$ 
ce qui prouve que $\phi_x$ est nul dans $\H^1(G, \Hom(\Omega, R))$. 
\qed

\begin{corollary}
\label{CalculDimensionF}
Soit $G$ un groupe $p$-cyclique agissant sur $X$ sans permuter les
branches et avec un conducteur $(m, m')$ ($m, m'\in \NN\cup
\{\infty\}$). Notons  
$$a=\begin{cases}
0 \textrm{ si } m < \infty \textrm{ et } m+1\not = 0 \mod p \\
1 \textrm{ sinon } \\
\end{cases}$$
$$a'=\begin{cases}
0 \textrm{ si } m' < \infty \textrm{ et } m'+1\not = 0 \mod p \\
1 \textrm{ sinon }. \\
\end{cases}$$
Alors le noyau du morphisme
$\Hom_R(\Omega, R'/R)^G  \to \H^1(G, \Hom_R(\Omega, R))$
est de dimension $a+a'.$
\end{corollary}

\preuve Le morphisme 
$\Hom_R(\Omega, R'/R)^G  \to \H^1(G, \Hom_R(\Omega, R))$ s'identifie
au morphisme composé 
$$k \d x \oplus k \d y \to k \phi_x \oplus k\phi_y \to \H^1(G,
\Hom_R(\Omega, R)).$$  
Le résultat est donc un corollaire immédiat de la proposition
précédente car $\phi_x$ est nul (en tant qu'application $G \to
\Hom(\Omega, R)$) si et seulement si $m=\infty$. \qed 

\begin{proposition}
\label{LocalSingulier2}
Soit $G$ un groupe cyclique agissant sur $X$ sans permuter les
branches et avec une différente $(\dif, \dif')$. Notons $G_{x, 0}$
l'image de $G$ dans $\Aut\, k((x))$. Alors l'image du morphisme 
$$\H^1(G, \Hom_{k[[x]]}(k[[x]]dx, k[[x]])) \to \H^1(G,
\Ext^0_R(\Omega, R'/R))$$ est nulle si $2 \dif +1=0 \mod |G_{x, 0}|$
et si $G_{x, 0} = G$, et de dimension 1 sinon. 
De plus, les images de
$\H^1(G, \Hom_{k[[x]]}(k[[x]]dx, k[[x]]))$ et $\H^1(G, \Hom_{k[[y]]}(k[[y]]dy, k[[y]]))$
dans $\H^1(G, \Ext^0_R(\Omega, R'/R))$ sont en somme directe. 
\end{proposition}

\preuve D'après le théorème \ref{ElementNonNulTraceNulle}, on peut
trouver un champ de vecteurs $f\d x \in \Hom_{k[[x]]}(k[[x]]dx,
k[[x]])$ avec $f(0) \not = 0$ et $\Tr_G\left(f \d x\right) =0$ si et
seulement si $2\dif+1\not = 0 \mod |G_{x, 0}|$ ou $G_{x, 0} \not =
G$. Supposons qu'il existe un tel élément. Comme $G$ agit trivialement
sur $\Ext^0_R(\Omega, R'/R)$, on voit que l'image de $f \d x$ est non
nulle dans $\H^1(G, \Ext^0_R(\Omega, R'/R))$.  

D'autre part, comme $G$ est cyclique et agit trivialement sur
$\Ext^0_R(\Omega, R'/R)$, on a $\H^1(G, \Ext^0_R(\Omega, R'/R)) \cong
\Ext^0_R(\Omega, R'/R)$, l'image est donc au plus de dimension $1$.  

\medskip

La dernière assertion est une conséquence immédiate de la description
de l'espace vectoriel $\Ext^0_R(\Omega, R'/R)$ provenant de l'isomorphisme
\eqref{IdentificationDualOmega} qui donne un isomorphisme canonique
$\H^1(G, \Ext^0_R(\Omega, R'/R))\cong k \d x \oplus k \d y.$
\qed

\begin{corollary}
\label{CalculCasCyclique}
 Soit $G$ un groupe $p$-cyclique agissant sur $X$ sans permuter les
branches et librement sur un ouvert dense (comme dans l'exemple
\ref{ActionpCycliqueSingulier}). Notons $(m, m')$ le conducteur de
$G$. Alors  
$$\dim_k \H^1(G, \Ext^0(\Omega, R))=m+m'+2+\left\lfloor {m \over p}
\right\rfloor-\left\lceil {2m+1 \over p} \right\rceil+\left\lfloor {m'
\over p} \right\rfloor-\left\lceil {2m'+1 \over p} \right\rceil.$$ 
\end{corollary}

\preuve C'est un conséquence immédiate du théorème
\ref{DimensionLisse} et des propositions \ref{LocalSingulier1} et
\ref{LocalSingulier2}. \qed 

\begin{remark}
On peut noter que, dans le cas général, le $k$-espace vectoriel
$\H^1(G, \Hom_{R'}(\Omega', R'))$ n'est pas de dimension finie (il
l'est si l'action de $G$ n'est triviale sur aucune branche). 
\end{remark}

\subsubsection{Épaississement de la singularité}
\label{EpaississementSingulariteLocal}

Nous allons maintenant nous intéresser à l'épaississement de la
singularité, c'est-à-dire à l'image du morphisme ${\textit{ép}}$ défini dans
la suite \eqref{EquationPointDouble}. Notre but est d'arriver au
théorème \ref{RelevementActionSingulier} qui décrit l'existence de
relèvements équivariants non topologiquement triviaux. 

Comme d'habitude, nous nous plaçons dans la cas où $G$ est un
$p$-groupe. On voit alors aisément que $\H^0(G,\Ext^1(\Omega, R))$ est
un $k$-espace vectoriel de dimension $1$ égal à $\Ext^{1}(\Omega, R)$.

Nous allons voir maintenant à quelles conditions un automorphisme peut 
se relever en un automorphisme d'une déformation non topologiquement 
triviale. Avant cela,  nous avons besoin de quelques lemmes 
préliminaires.  

\begin{lemma}
\label{RelevementMorphisme}
Soient $\sigma$ un automorphisme de $X$ ne permutant pas les bran\-ches.  Notons 
$$\X=\Spec k[\epsilon]/(\epsilon^2)[[x, y]]/(xy-\epsilon).$$ 
Soit $\sigma_\epsilon$ un relèvement de $\sigma$ à $\X$, alors il
existe $h_{\sigma, 0} \in xk[[x]]$ et $h_{\sigma, 1}\in yk[[y]]$ tels
que 
\begin{eqnarray}
\label{RelevementX}
\sigma_\epsilon(x)=\sigma(x)-\epsilon {\sigma(y)- y \over y\sigma(y)}+\epsilon h_{\sigma, 0}(x), \\
\label{RelevementY}
\sigma_\epsilon(y)=\sigma(y)-\epsilon {\sigma(x)- x \over x\sigma(x)}+\epsilon h_{\sigma, 1}(y).
\end{eqnarray}
De plus, pour tout $h_{\sigma, 0} \in xk[[x]]$ et $h_{\sigma, 1} \in ta[[y]]$, il existe un automorphisme $\sigma_\varepsilon:\X \rightarrow \X$ donné par les équations \eqref{RelevementX} et \eqref{RelevementY}.
\end{lemma}

\preuve Notons encore $\sigma(x)$ et $\sigma(y)$ les relèvements
triviaux de $\sigma(x)$ et $\sigma(y)$ à $k[\epsilon][[x]]$ et
$k[\epsilon][[y]]$. Comme $\sigma_\epsilon$ relève $\sigma$, on a
$\sigma_\epsilon(x)-\sigma(x)=0 \mod \varepsilon$. Par suite, il
existe $h_{\sigma, 0} \in xk[[x]]$ et $e_{\sigma, 0} \in k[[y]]$
uniques tels que  
$$\sigma_\epsilon(x)-\sigma(x)=\epsilon(e_{\sigma, 0}(y)-h_{\sigma, 0}(x)).$$ De même, il existe $h_{\sigma, 1} \in yk[[y]]$ et $e_{\sigma, 1} \in k[[x]]$ uniques tels que 
$$\sigma_\epsilon(y)-\sigma(y)=\epsilon(e_{\sigma, 1}(x)-h_{\sigma, 1}(y)).$$
Or $\sigma_\epsilon$ est un automorphisme de $k[\epsilon][[x,
y]]/(xy-\epsilon)$ donc $\sigma_\epsilon(xy-\epsilon)=0 \mod
(xy-\epsilon)$. En exprimant cette égalité, on obtient les équations
\eqref{RelevementX} et \eqref{RelevementY} en faisant un développement
limité à l'ordre $1$.
La réciproque se démontre de la même manière. \qed

Dans le cas de l'action d'un groupe, les séries formelles $h_{\sigma,
0}$ et $h_{\sigma, 1}$ du lemme précédent doivent vérifier certaines
conditions de compatibilité lorsque $\sigma$ varie dans $G$. Pour des
raisons techniques, il est plus commode de poser $f_{\sigma,
0}(x)=(\der x \sigma(x)) h_{\sigma, 0}(x)$ et $f_{\sigma, 1}(y)=(\der
y \sigma(y)) h_{\sigma, 1}(y)$. 

\begin{proposition}
\label{RelevementActionGroupe}
Soit $G$ un $p$-groupe fini agissant sur $X$ sans permuter les
bran\-ches. Pour tout $\sigma \in G$,
notons $v_{\sigma, 0}$ (resp. $v_{\sigma, 1}$) le terme de degré $2$
dans la décomposition de $\sigma(x)$ (resp. $\sigma(y)$). Alors il
existe un relèvement de l'action de $G$ à $\X$ si et seulement s'il
existe des séries formelles $f_{\sigma, 0} \in xk[[x]]$ et $f_{\sigma,
1} \in yk[[y]]$ pour tout $\sigma \in G$ vérifiant 
\begin{eqnarray}
\label{RecurrenceX}
f_{\sigma\tau, 0}(x)=\left(1-{1 \over ( \der x \tau)(\sigma(x))}\right){v_{\sigma, 1} \over\der x \sigma(x)}+f_{\sigma, 0}(x) +{1 \over \der x \sigma(x)}f_{\tau, 0}(\sigma(x)) \\
\label{RecurrenceY}
f_{\sigma\tau, 1}(y)=\left(1-{1\over ( \der y \tau)(\sigma(y))}\right){v_{\sigma, 0} \over \der y \sigma(y)}+f_{\sigma, 1}(y)+{1 \over \der y \sigma(y)}f_{\tau, 1}(\sigma(y))
\end{eqnarray}
où $(\der x \tau)(\sigma(x))$ désigne l'évaluation de la série formelle $\der x (\tau(x))$ en $\sigma(x)$.
\end{proposition}

\preuve Montrons que c'est une condition nécessaire. Supposons qu'il
existe un relèvement de l'action de $G$ à $\X$. Pour tout $\sigma\in
G$ notons $\sigma_\epsilon$ le relèvement de $\sigma$ à $\X$ défini
par le relèvement de $G$. D'après le lemme \ref{RelevementMorphisme},
pour tout $\sigma \in G$ il existe $h_{\sigma, 0} \in xk[[x]]$ et
$h_{\sigma, 1} \in yk[[y]]$ telles que le relèvement $\sigma_\epsilon$
de $\sigma$ à $\X$ vérifie les équations \eqref{RelevementX} et
\eqref{RelevementY}. 

De plus, comme on a supposé que $G$ agit sur $\X$, on a, pour tout $\sigma, \tau \in G$
$$\left(\sigma \circ \tau\right )_\epsilon=\sigma_\epsilon \circ \tau_\epsilon.$$
Un développement limité à l'ordre 1 (direct mais assez pénible à cause des compositions de morphismes) permet de trouver
$$h_{\sigma\tau, 0}(x)=((\der x \tau)(\sigma(x))-1) {v_{\sigma, 1}}+h_{\sigma, 0}(x)(\der x \tau)(\sigma(x))+h_{\sigma, 0}(\sigma(x)). $$
Par suite, en posant $f_{\rho, 0}(x)={h_{\rho, 0}(x) \over \der x \rho(x)}$ pour tout $\rho \in G$ on trouve la formule \eqref{RecurrenceX}. Des calculs similaires conduisent à formule \eqref{RecurrenceY}.
Ces mêmes calculs montrent que les conditions \eqref{RecurrenceX} et \eqref{RecurrenceY} sont suffisantes pour relever l'action. \qed

\begin{proposition}
\label{RelevementActionGroupeCyclique}
Reprenons les notations de la proposition précédente et supposons de
plus que $G$ est cyclique d'ordre $n$. Soit $\sigma$ un générateur de
$G$. Alors il existe un relèvement de l'action de $G$ à $\X$ si et
seulement s'il existe $f_{\sigma, 0} \in xk[[x]]$ et $f_{\sigma, 1}
\in yk[[y]]$ telles que 
\begin{equation}
\label{FormuleTrace}
\Tr_G \left(\left( f_{\sigma, 0}(x) +v_{\sigma, 1}\right)\d x \right )= 0= 
\Tr_G \left(\left(f_{\sigma, 1}(y) +v_{\sigma, 0} \right)\d y \right ).
\end{equation}
Si de tels éléments existent alors le relèvement est donné par les
formules \eqref{RelevementX} et \eqref{RelevementY} avec $h_{\sigma,
0}(x)={f_{\sigma, 0}(x) \der x \sigma(x)}$ et $h_{\sigma,
1}(y)={f_{\sigma, 1}(y) \der y \sigma(y) }$. 
\end{proposition}

\preuve Supposons qu'on puisse relever l'action de $G$ à $\X$. Soit $i
\in \NN$, on a, d'après la formule \eqref{RecurrenceX}, 
$$f_{\sigma^{i+1}, 0}(x)=\left(1-{1 \over ( \der x
\sigma^i)(\sigma(x))}\right){v_{\sigma, 1} \over \der x
\sigma(x)}+f_{\sigma, 0}(x) +{1 \over \der x \sigma(x)}f_{\sigma^i,
0}(\sigma(x)).$$ 
Une récurrence immédiate montre que
$$f_{\sigma^{i+1}, 0}(x) = v_{\sigma, 1} \left(\sum_{j=0}^i {1 \over
\der x \sigma^{j+1}(x)}-{i \over \der x \sigma^{i+1}(x)}\right
)+\sum_{j=0}^i {f_{\sigma, 0}(\sigma^j(x)) \over \der x
\sigma^{j}(x)}.$$ 
C'est-à-dire
\begin{multline*}
f_{\sigma^{i+1}, 0}(x)\d x = v_{\sigma, 1}\sum_{j=0}^i
\sigma^{j+1}\left (\d x \right )- {i \over \der x \sigma^{i+1}(x)}{\d
x}+\sum_{j=0}^i \sigma^j\left(f_{\sigma, 0}(x)\d x\right). 
\end{multline*}
Or par définition de $n$ on a $\sigma^n=Id$ donc $f_{\sigma^n, 0}=0$ ce qui impose 
$$0 = v_{\sigma, 1}\sum_{j=0}^{n-1} \sigma^{j+1}\left (\d x \right
)+\sum_{j=0}^{n-1} \sigma^j\left(f_{\sigma, 0}(x)\d x\right)$$ 
ce qui est bien le résultat annoncé. Des calculs similaires conduisent
à l'équation sur l'autre branche. 

\medskip
Réciproquement, supposons qu'il existe $f_{\sigma, 0}$ et $f_{\sigma,
1}$ vérifiant les hypothèses de la proposition. Posons, pour tout $i
\in \{0, \ldots, n-1\}$,  
$$f_{\sigma^{i+1}, 0}(x):= v_{\sigma, 1} \left(\sum_{j=0}^i {1 \over
\der x \sigma^{j+1}(x)}-{i \over \der x \sigma^{i+1}(x)} \right
)+\sum_{j=0}^i {f_{\sigma, 0}(\sigma^j(x)) \over \der x \sigma^{j}(x)}$$ 
(et une formule semblable pour $f_{\sigma^{i+1}, 1}$). Par
construction, les $f_{\sigma^i, 0}$ et $f_{\sigma^i, 1}$ vérifient les
équations \eqref{RecurrenceX} et \eqref{RecurrenceY}. La proposition
\ref{RelevementActionGroupe} permet alors de relever l'action de $G$ à
$\X$. \qed 

Nous aurons besoin, lors de l'étude globale, d'une description de
l'action sur chaque branche d'une déformation du premier ordre du
point double. Ceci est donné dans le lemme suivant. 

\begin{lemma}
\label{ActionSurAnneauDeFraction}
Soit $\X \to \Spec k[\epsilon]/(\epsilon^2)$ une déformation
$G$-équivariante de $X$. Re\-prenons les hypothèses et notations de la
proposition \ref{RelevementActionGroupe}. Alors pour tout $\sigma \in
G$ il existe $h_{\sigma, 0} \in xk[[x]]$, $h_{\sigma, 1} \in yk[[y]]$
et $\lambda \in k$ tels que l'action sur $k[\epsilon]((x))$ et sur
$k[\epsilon]((y))$ soit donnée par 
$$\sigma_\epsilon(x)=\sigma(x)-\lambda \epsilon v_{\sigma, 1}+\epsilon h_{\sigma, 0}(x)$$
$$\sigma_\epsilon(y)=\sigma(y)-\lambda \epsilon v_{\sigma, 0}+\epsilon h_{\sigma, 1}(y).$$
De plus, $\lambda = 0$ si et seulement si $\X$ est une déformation
topologiquement triviale. 
\end{lemma}

\preuve C'est une conséquence immédiate du lemme \ref{RelevementMorphisme}. \qed 

En particulier, on peut constater que les séries de Laurent
$\sigma_\epsilon(x)$ et $\sigma_\epsilon(y)$ n'ont pas de pôles. Cette
remarque sera très utile lors de l'étude des déformations globales. 

\begin{definition}
\label{ActionRelevable}
Soit $G$ un groupe agissant sur $X=\Spec k[[x, y]]/(xy)$. On dira que
l'action de $G$ est \textit{relevable à l'ordre 1} s'il existe une
déformation $G$-équivariante de $X$ non topologiquement triviale. Dans
le cas contraire, elle sera dites \textit{non relevable à l'ordre 1}.  

Supposons que l'action est relevable à l'ordre 1. Notons $D_x$
(resp. $D_y$) le stabilisateur de $(x)$ (resp. $(y)$), $G_x$
(resp. $G_y$) l'image dans $\Aut_k(\Spec k((x)))$ (resp. $\Aut_k(\Spec
k((y)))$) de $D_x$ (resp. $D_y$) et $T_x$ (resp. $T_y$) le noyau du
morphisme $D_x \to G_x$ (resp. $D_y \to G_y$). Nous dirons que
l'action de $G$ est \textit{relevable inconditionnellement} s'il
existe une déformation $G$-équivariante $\X$ de $X$ non
topologiquement triviale telle que l'action induite de $T_x$ sur
$\Spec k[\epsilon]((x))$ (resp. $T_y$ sur $\Spec k[\epsilon]((y))$)
est triviale. Dans le cas contraire, l'action sera dite
\textit{relevable conditionnellement}. 
\end{definition}

\begin{theorem}
\label{RelevementActionSingulier}
Soit $G$ un $p$-groupe agissant sur $X=\Spec k[[x,y]]/(xy)$ sans
permuter les branches. Reprenons les notations de la définition
\ref{ActionRelevable}. 
\begin{enumerate}[{\rm(i)}]
\item Si le conducteur est différent de $(1, m)$ et $(m, 1)$ (avec $m
\in \NN\cup\{\infty\}$), alors l'action de $G$ est relevable inconditionnellement. 
\item Supposons que $G$ est un $p$-groupe cyclique et que son
conducteur est $(m, 1)$ (avec $m \in \NN \cup 
\{\infty\}$ et $m \not = 1$). Notons $(\dif, \dif')$ sa
différente. Alors l'action de $G$ est relevable si et seulement si
$2\dif+1 \not = 0 \mod |G_x|$ ou $G_x \not = G$. 
De plus, si l'action est relevable, alors elle l'est
inconditionnellement si et seulement si $2\dif+1 \not = 0 \mod |G_x|$.  
\item Si $G$ est un $p$-groupe cyclique agissant sur $X$ avec un
conducteur $(1, 1)$, alors l'action est 
relevable si et seulement si  
$$\begin{cases}
2\dif+1 \not = 0 \mod |G_x|$ ou $G_x \not = G \\
2\dif'+1 \not = 0 \mod |G_y|$ ou $G_y \not = G.
\end{cases}$$
De plus, l'action est relevable inconditionnellement si et seulement
si $2\dif+1 \not = 0 \mod |G_x|$ et $2\dif'+1 \not = 0 \mod |G_y|$. 
\end{enumerate}
\end{theorem}

\preuve Montrons d'abord le point (i). Pour tout $\sigma \in G$ on
voit que $v_{\sigma, 0} = 0$ et $v_{\sigma, 1}=0$. Par suite, posons
$f_{\sigma,0}=f_{\sigma, 1}=0$ pour tout $\sigma \in G$. On voit alors
que $\{f_{\sigma, 0}\}_{\sigma \in G}$ est solution des équations
\eqref{RecurrenceX} et $\{f_{\sigma, 1}\}_{\sigma \in G}$ des
équations \eqref{RecurrenceY} et donc qu'il existe un relèvement
inconditionnel (par le lemme \ref{ActionSurAnneauDeFraction}) de $X$. 

\medskip

Passons à la démonstration du point (ii). Supposons dans un premier
temps que $2\dif+1 \not = 0 \mod |G_x|$ ou $G_x \not = G$. 
Soit $f_{0} \in k[[x]]$ un élément de trace nulle fourni par le
théorème (\ref{ElementNonNulTraceNulle}). Alors la proposition
\ref{RelevementActionGroupeCyclique} fournit un relèvement (avec
$f_{\sigma, 0}=f_0$ et $f_{\sigma, 1}=0$ où $\sigma$ est un générateur
de $G$).  
D'autre part, la proposition \ref{RelevementActionGroupeCyclique} et
le théorème \ref{ElementNonNulTraceNulle} fournissant des
équivalences, la réciproque est démontrée. 

\medskip

Revenons au cas $2\dif+1 \not = 0 \mod |G_x|$. Alors on peut supposer
que $\Tr_{G_x} f_0 = 0$ d'après le théorème
\ref{ElementNonNulTraceNulle}. On voit ainsi que l'action précédemment
définie induit une action de $G_x$ sur $\Spec k[\epsilon]((x))$ et
donc que l'action est relevable inconditionnellement. 

Réciproquement, supposons que l'action est relevable
inconditionnellement. Soit $(\X, G)$ une déformation équivariante de
$(X, G)$ non topologiquement triviale et telle que $G_x$ agisse sur
$\Spec k[\epsilon]((x))$. Notons $f_0$ l'élément donné par la
proposition \ref{RelevementActionGroupeCyclique}. Comme $G_x$ agit sur
$\Spec k[\epsilon]((x))$, on montre aisément que $\Tr_{G_x} f_0=0$ et
$f_0(0) \not = 0$ (comme dans la proposition
\ref{RelevementActionGroupeCyclique}). Le résultat vient alors du
théorème \ref{ElementNonNulTraceNulle}. 

\medskip

La démonstration du point (iii) est similaire à celle du point (ii).\qed

\section{Étude des déformations globales}

Cette section a pour but l'étude des déformations équivariantes des
courbes semi-stables. Nous ne supposons pas, excepté dans la
sous-section \ref{DeformationCourbeLisse}, que l'action est libre sur
un ouvert dense. 

\subsection{Déformation des courbes lisses}
\label{DeformationCourbeLisse}
Soient $k$ un corps algébriquement clos, $C \to \Spec k$ une courbe
lisse et $G$ un groupe agissant fidèlement sur $C$ (en particulier,
l'action est libre sur un ouvert dense). Notons $C_{ram}$ l'ensemble
des points possédant un stabilisateur non trivial. Pour tout $\p \in
C_{ram}$ nous noterons $\dif_\p$ la différente de $G$ en $\p$. Pour
chaque classe $\alpha \in C_{ram}/G$, choisissons un
représentant $\p_\alpha$ dans la classe $\alpha$. On a une suite
exacte provenant de la suite spectrale $\textrm{I}_{2}^{p,
q}:=\H^p(C/G, \ext_{G}^{q} (\Omega_{X/k}, \O_C))$
\begin{multline*}
0 \to \H^1(C/G, \ext^0_G(\Omega_{C/k}, \O_C)) \to \Ext^1_G(\Omega_{C/k}, \O_C) \to \\
\bigoplus_{\alpha \in C_{ram}/G}
\Ext^1_{D_{\p_\alpha}}(\hat\Omega_{C/k, \p_{\alpha}}, \hat\O_{C,
\p_\alpha}) \to 0 
\end{multline*}
(cf. \cite{BertinMezard}, Lemme 3.3.2). Un calcul similaire à celui de
la proposition 5.3.2 de \cite{BertinMezard} permet alors de montrer
que 
$$\dim_k \H^1(C/G, \ext^0_G(\Omega_{C/k}, \O_C)) =
3p_a(C/G)-3+\sum_{\alpha \in C_{ram}/G} \left\lceil { \dif_{\p_\alpha}
\over |D_{\p_\alpha}|} \right \rceil$$ 
(on prendra garde à la coquille qui s'est glissée dans l'énoncé et
dans la preuve de la proposition 5.3.2 de \cite{BertinMezard} et qui
consiste en la transformation des parties entières supérieures en
parties entières inférieures). 
Dans le cas où le stabilisateur en chaque point est cyclique on peut
donc calculer, grâce à \cite{BertinMezard}, Théorème 4.1.1, la
dimension de $\Ext^1_G(\Omega_{C/k}, \O_C)$ explicitement, ce qui
donne une généralisation immédiate de la proposition 5.3.2 de
\cite{BertinMezard}. 
\begin{theorem}
\label{CalculGlobalCasLisse}
Soient $C\to \Spec k$ une courbe lisse sur un corps algébriquement
clos de caractéristique $p$ et $G$ un $p$-groupe fini agissant
fidèlement sur $C$ et tel que pour tout point fermé $\p \in C$ le
stabilisateur $D_\p$ soit cyclique. Alors 
$$\dim_k \Ext^1_G(\Omega_{C/k}, \O_C)=3p_a(C/G)-3+\sum_{\alpha \in
C_{ram}/G} \left\lfloor { 2\dif_{\p_\alpha} \over |D_{\p_\alpha}|}
\right \rfloor.$$ 
\end{theorem}

\preuve D'après la suite exacte précédente, on a
\begin{multline*}
\dim_k \Ext^1_G(\Omega_{C/k}, \O_C)=\dim_k \H^1(C/G,
\ext^0_G(\Omega_{C/k}, \O_C))+ \\\sum_{\alpha \in C_{ram}/G} \dim_k
\Ext^1_{D_{\p_\alpha}}(\hat\Omega_{C/k, \p_{\alpha}}, \hat\O_{C,
\p_\alpha}). 
\end{multline*}
Finalement, grâce au théorème \ref{DimensionLisse}, on obtient le résultat. \qed

\subsection{Déformation des courbes singulières}
\label{DeformationCourbeSinguliere}

Soient $C \to \Spec k$ une courbe semi-stable et $G$ un groupe fini
d'automorphismes de $C$.  
La suite spectrale $\textrm{I}_{2}^{p,q}:=\H^p(C/G, \ext_{G}^{q}
(\Omega_{X/k}, \O_C))$ donne une suite exacte 
\begin{multline*}
0 \to \H^1(C/G, \ext^0_G(\Omega_{C/k}, \O_C)) \to \Ext^1_G(\Omega_{C/k}, \O_C) \to  \\ \H^0(C/G, \ext^1_G(\Omega_{C/k}, \O_C)) \to  \H^2(C/G, \ext^0_G(\Omega_{C/k}, \O_C)).
\end{multline*}
Comme $C/G$ est une courbe, on a $\H^2(C/G, \ext^0_G(\Omega_{C/k}, \O_C))=0$. Par suite, calculer la dimension de $\Ext^1_G(\Omega_{C/k}, \O_C)$ revient à calculer celles des deux autres termes de la suite.

Considérons le diagramme commutatif suivant
$$\xymatrix{
\tilde C \ar[r]^\varphi\ar[d]_{\tilde \pi} & C \ar[d]^\pi \\
\tilde C/G \ar[r]_\psi & C/G,
}$$ les morphismes $\varphi$ et $\psi$ étant les morphismes de normalisation et les morphismes $\pi$ et $\tilde \pi$ étant les morphismes quotients.

La suite spectrale $\textrm{II}_{2}^{p, q}:=\R^p\pi_{*}^{G}\ext^{q}
(\Omega_{C/k}, \O_C)$ donne  une suite exacte
$$0 \to \left(\R^1 \pi_*^G \right)\Omega_{C/k}^\vee \to \ext^1_G(\Omega_C, \O_C) \to \pi_*^G(\ext^1(\Omega_{C/k}, \O_C)).$$
On a donc une suite exacte
\begin{multline*}
0 \to \H^0(C/G, \left(\R^1 \pi_*^G \right)\Omega_{C/k}^\vee) \to \H^0(C/G, \ext^1_G(\Omega_C, \O_C)) \to \\
\H^0(C, \ext^1(\Omega_{C/k}, \O_C))^G.
\end{multline*}
En résumé, on a le diagramme commutatif suivant (à ligne et colonne exactes)
\begin{footnotesize}
$$\xymatrix@C0.3cm@R0.5cm{
 & & & 0 \ar[d] & \\
 & & & \H^0(C/G, \R^1 \pi_*^G \Omega_{C/k}^\vee)\ar[d] & \\
0 \ar[r] & \H^1(C/G, \ext^0_G(\Omega_{C/k}, \O_C)) \ar[r] & \Ext^1_G(\Omega_{C/k}, \O_C) \ar[r]\ar[rd]_{\textit{ép}} & \H^0(C/G, \ext^1_G(\Omega_{C/k}, \O_C)) \ar[r]\ar[d] & 0 \\
 & & & \H^0(C/G, \pi_*^G \ext^1(\Omega_{C/k}, \O_C)) & 
}$$
\end{footnotesize}

Nous allons relier 
$\H^1(C/G, \ext^0_G(\Omega_{C/k}, \O_C))$ et $\H^0(C/G, \R^1
\pi_*^G\Omega_{C/k}^\vee)$ à des groupes similaires définis à partir
de $\tilde C$ (comme dans le cas local). Nous considérerons ensuite l'image du morphisme  
$${\textit{ép}}:\Ext^1_G(\Omega_{C/k}, \O_C) \to \H^0(C, \ext^1(\Omega_{C/k},
\O_C))^G.$$ 

\subsubsection{Lien avec la normalisée}
\label{LienNormalisee}

Le but de cette partie est de montrer le théorème \ref{LienNormaliseGlobale}. Avant cela, nous avons besoin d'introduire quelques notations.

Pour tout schéma $Z$ nous noterons $\pi_0(Z)$ l'ensemble des composantes connexes de $Z$.
Soient $C \to \Spec k$ une courbe stable et $G$ un $p$-groupe fini
d'automorphismes de $C$. Pour tout $\beta \in \pi_0(\tilde C)/G$,
choisissons un représentant $C_\beta$ dans la classe $\beta$ et notons
$D_\beta$ son stabilisateur, $T_\beta$ le noyau de l'homo\-morphisme  
$D_\beta \to \Aut(\Spec \kappa(C_\beta))$ et
$G_\beta:=D_\beta/T_\beta$. Notons de plus $\varphi_\beta:C_\beta \to
C$ le morphisme induit par la normalisation, $\pi_\beta:C_\beta \to
C_\beta / G_\beta$ le morphisme quotient et
$\psi_\beta:C_\beta/G_\beta \to C/G$. 
On a donc un diagramme commutatif
\begin{equation}\label{DefinitionMorphismesGlobaux}
\xymatrix{
C_\beta \ar@/^0.5cm/[rr]^{\varphi_\beta}\ar[r]\ar[d]_{\pi_\beta} & \tilde C \ar[r]_\varphi\ar[d]_{\tilde \pi} & C\ar[d]_{\pi} \\
C_\beta/G_\beta \ar[r]\ar@/_0.5cm/[rr]_{\psi_\beta} & \tilde C/G \ar[r]^\psi & C/G
}
\end{equation}

La suite exacte \eqref{MorphismeNormaliseeLocal} possède un analogue global (dont l'exactitude se démontre de manière identique)
\begin{equation}
\label{SuiteLienNormaliseGlobal}
0 \to \Omega_{C/k}^\vee \to \varphi_* \Omega_{\tilde C/k}^\vee \to \hom(\Omega_{C/k}, \left( \varphi_* \O_{\tilde C} \right ) / \O_C) \to 0.
\end{equation}
Notons $$\Gg_\beta:=\Im\left(\R^1 \pi_{\beta, *}^{G_\beta}
\left(\varphi_\beta^* \varphi_* \Omega_{\tilde C/k}^\vee \right) \to
\R^1 \pi_{\beta, *}^{G_\beta} \left(\varphi_\beta^*
\hom(\Omega_{\tilde C/k}, \varphi_* \O_{\tilde C}/\O_C)
\right)\right)$$
(on remarque que les groupes
apparaissant dans la définition de $\Gg_\beta$ sont les $G_\beta$ et
non les $D_\beta$). 
On voit que le support de $\Gg_\beta$ est dans $\psi_\beta^{-1}(\pi(C_{sing}))$.

\begin{exemple}
\label{CalculExpliciteG}
Supposons que pour tout $\p \in C_{sing}$, $D_\p$ est cyclique et ne
permute pas les branches. Pour tout $\alpha \in C_{sing}/G$
choisissons un point $\p_\alpha$ dans la classe de $\alpha$ et notons
$(\dif_\alpha, \dif'_\alpha)$ la différente de $D_{\p_\alpha}$. La
proposition \ref{LocalSingulier2} permet de voir que  
$$\sum_{\beta} \dim_k \H^0(C_\beta/G_\beta,
\Gg_\beta)=\left|\left\{\substack{\alpha \in C_{sing}/G \\
2\dif_\alpha+1 \not = 0 \mod |G_{\p_\alpha}|\\ \textrm{ou }
G_{\p_\alpha} \not =
D_{\p_\alpha}}\right\}\right|+\left|\left\{\substack{\alpha \in
C_{sing}/G \\ 2\dif'_\alpha+1 \not = 0 \mod |G_{\p_\alpha}|\\
\textrm{ou } G_{\p_\alpha} \not = D_{\p_\alpha}}\right\}\right|.$$  
\end{exemple}

De manière analogue, on définit
$$\Fg:= \ker \left(\pi_*^G\hom(\Omega_{C/k}, \left( \varphi_* \O_{\tilde C} \right ) / \O_C) \to \R^1 \pi_*^G \Omega_{C/k}^\vee \right).$$ \index{${\mathfrak F}$}
On voit que $\Fg$ est à support dans $\pi(C_{sing})$.
On peut préciser la structure des fibres de $\Fg$ dans le cas général. Ceci est donné par la proposition suivante.

\begin{proposition}
Reprenons les notations précédentes. Soient $\q$ un point de $C/G$ et
$\p \in \pi^{-1}(\q)$.  Notons $D_\p$ le stabilisateur de $\p$, $T_\p$
le noyau de l'homomorphisme $D_\p \to \Aut(\Frac(\O_{C, \p}))$ et
$G_\p:=D_\p /T_p$. Alors  
$$\Fg_\q = \ker \left(\hom(\Omega_{C/k}, \varphi_* \O_{\tilde
C}/\O_C)_\p^{G_\p} \to  \H^1(G_\p, \Omega_{C/k, \p}^\vee)\right).$$ 
\end{proposition}

\preuve La définition de $\Fg$ et une décomposition selon les orbites
sous $D_{\p}$ (cf. \cite{CohomologyOfGroup}, 27.1)  permettent de voir que 
$$\Fg_\q = \ker \left(\hom(\Omega_{C/k}, \varphi_* \O_{\tilde
C}/\O_C)_\p^{G_\p} \to  \H^1(D_\p, \Omega_{C/k, \p}^\vee)\right).$$ 
Il est alors aisé de voir que le morphisme se factorise en fait par
$\H^1(G_\p, \Omega_{C/k, \p}^\vee)$ car $T_\p$ agit trivialement.
 \qed

\begin{exemple}
\label{CalculExpliciteF}
Supposons que pour tout $\p \in C_{sing}$ le stabilisateur $D_\p$ est
$p$-cyclique et agit sans échanger les branches. Pour tout $\alpha \in
C_{sing}/G$ choisissons un point $\p_\alpha$ dans la classe $\alpha$
et notons $(m_\alpha, m'_\alpha)$ le conducteur de
$D_{\p_\alpha}$. Le corollaire \ref{CalculDimensionF} nous donne 
$$\dim_k \H^0(C, \Fg)=\left|\left\{\substack{\alpha \in C_{sing}/G \\
m_\alpha+1 = 0 \mod p \\ \textrm{ou } m_\alpha =
\infty}\right\}\right|+\left|\left\{\substack{\alpha \in C_{sing}/G \\
m'_\alpha+1 = 0 \mod p \\ \textrm{ou } m'_\alpha =
\infty}\right\}\right|.$$ 
\end{exemple}

\begin{theorem}
\label{LienNormaliseGlobale}
Soient $C \to \Spec k$ une courbe stable et $G$ un groupe fini agissant sur $C \to \Spec k$.
Reprenons les notations précédentes.
Alors on a deux suites exactes
\begin{multline}
\label{SuiteNormaliseeTermeGlobal}
0 \to \Ext^0_G(\Omega_{\tilde C/k}, \O_{\tilde C}) \to \H^0(C/G, \Fg)
\to \H^1(C/G, \ext^0_G(\Omega_{C/k}, \O_C)) \to  \\ \H^1(\tilde C/G,
\ext^0_G(\Omega_{\tilde C/k}, \O_{\tilde C})) \to 0,  
\end{multline}
\begin{multline}
\label{SuiteNormaliseeTermeLocal}
\hskip -0.3cm 0 \to \H^0(C/G, \Fg) \to \H^0(C, \hom(\Omega_{C/k},
\left( \varphi_* \O_{\tilde C} \right ) / \O_C))^G \to  \H^0(C/G, \R^1
\pi_*^G \Omega_{C/k}^\vee) \\  \to \bigoplus_{\beta\in\pi_0(\tilde
C)/G}\H^0(C_\beta/G_\beta, \R^1 \tilde\pi_*^{G_\beta}
\Omega_{C_\beta/k}^\vee) \to \bigoplus_{\beta\in\pi_0(\tilde C)/G}
\H^0(C_\beta/G_\beta, \Gg_\beta) \to 0. 
\end{multline}
\end{theorem}

\preuve
Nous allons mener une étude parallèle à celle de la section
\ref{LienNormaliseeLocal}. 

\medskip

En appliquant le foncteur $\pi_*^G$ à la suite
\eqref{SuiteLienNormaliseGlobal}, on obtient une suite exacte longue  
\begin{multline}
\label{SuiteExacteLongueGroupeFaisceaux}
0 \to \pi_*^G \Omega_{C/k}^\vee \to \psi_* \circ \tilde \pi_*^G
\Omega_{\tilde C/k}^\vee \to \pi_*^G\hom(\Omega_{C/k}, \left(
\varphi_* \O_{\tilde C} \right ) /\O_C) \to \left ( \R^1 \pi_*^G
\right )\Omega_{C/k}^\vee \\  
\to \psi_* \circ \left ( \R^1 \tilde \pi_*^G \right )\Omega_{\tilde
C/k}^\vee \to \left ( \R^1 \pi_*^G \right )\hom(\Omega_{C/k}, \left(
\varphi_* \O_{\tilde C} \right ) / \O_C) 
\end{multline}
car pour tout $i \ge 0$ on a
$$\psi_* \circ \R^i \tilde \pi_*^G=
\R^i\left(\psi \circ \tilde \pi \right)^G_*=
\R^i\left(\pi \circ \phi \right)_*^G =
(\R^i \pi_*^G) \circ \phi_*
$$
(ce qui découle de l'exactitude de $\varphi_*$ et de $\psi_*$).
Notons $\Hc$ le noyau du morphisme
$$\psi_* \circ \left ( \R^1 \tilde \pi_*^G \right )\Omega_{\tilde
C/k}^\vee \to \left ( \R^1 \pi_*^G \right )\hom(\Omega_{C/k}, \left(
\varphi_* \O_{\tilde C} \right ) / \O_C)$$  
de sorte qu'on a trois suites exactes
\begin{equation}
\label{SuiteExacteF}
0 \to \pi_*^G \Omega_{C/k}^\vee\to \psi_* \circ \tilde \pi_*^G
\Omega_{\tilde C/k}^\vee \to \Fg \to 0,  
\end{equation}
\begin{equation}
\label{SuiteExacteFH}
0 \to \pi_*^G\hom(\Omega_{C/k}, \left( \varphi_* \O_{\tilde C} \right ) / \O_C)/\Fg \to \left ( \R^1 \pi_*^G \right )\Omega_{C/k}^\vee \to \Hc \to 0,
\end{equation}
\begin{equation}
\label{SuiteExacteHG}
0 \to \Hc \to 
\psi_* \circ \left ( \R^1 \tilde \pi_*^G \right )\Omega_{\tilde
C/k}^\vee \to \ \left ( \R^1 \pi_*^G \right )\hom(\Omega_{C/k}, \left(
\varphi_* \O_{\tilde C} \right ) / \O_C). 
\end{equation}

En prenant la suite exacte longue de cohomologie des faisceaux
associée à la suite \eqref{SuiteExacteF} on obtient une suite exacte 
\begin{multline}
\label{SuiteExacteFGlobale}
0 \to \Ext^0_G(\Omega_{C/k}, \O_C) \to \Ext^0_G(\Omega_{\tilde C/k},
\O_{\tilde C}) \to \H^0(C/G, \Fg) \to  \\ \H^1(C/G,
\ext^0_G(\Omega_{C/k}, \O_C)) \to \H^1(\tilde C/G,
\ext^0_G(\Omega_{\tilde C/k}, \O_{\tilde C})) \to 0  
\end{multline}
car, $\Fg$ étant à support de dimension $0$, on a $$\H^1(C/G, \Fg)=0.$$
Par suite, on trouve la suite \eqref{SuiteNormaliseeTermeGlobal} car
$C$ est stable et donc $\Ext^0(\Omega_{C/k}, \O_C)=0$ d'après
\cite{Deligne_Mumford}, Lemma 1.4. 

D'autre part, la suite exacte longue de cohomologie associée à la
suite \eqref{SuiteExacteFH} nous donne, comme $\Fg$ et
$\pi_*^G\hom(\Omega_{C/k}, \left( \varphi_* \O_{\tilde C} \right ) /
\O_C)$ sont à support de dimension $0$, une suite exacte 
\begin{multline*}
0 \to \H^0(C/G, \Fg) \to \H^0(C, \hom(\Omega_{C/k}, \left( \varphi_*
\O_{\tilde C} \right ) / \O_C))^G \to  \\ \H^0(C/G, \R^1 \pi_*^G
\Omega_{C/k}^\vee) \to \H^0(C/G, \Hc) \to 0. 
\end{multline*}
De même, la suite exacte \eqref{SuiteExacteHG} nous donne une suite
exacte 
\begin{multline*}
0 \to \H^0(C/G, \Hc) \to \H^0(\tilde C/G, \R^1 \tilde\pi_*^G
\Omega_{\tilde C/k}^\vee) \to  \\ \H^0\left(C/G, \left ( \R^1 \pi_*^G
\right )\hom(\Omega_{C/k}, \left( \varphi_* \O_{\tilde C} \right ) /
\O_C)\right). 
\end{multline*}
On a donc une suite exacte
\begin{multline}
\label{SuiteNormalisee2}
\hskip-0.4cm 0 \to \H^0(C/G, \Fg) \to \H^0(C, \hom(\Omega_{C/k}, \left( \varphi_*
\O_{\tilde C} \right ) / \O_C))^G \to  \H^0(C/G, \R^1 \pi_*^G
\Omega_{C/k}^\vee)\\ \to \H^0(\tilde C/G, \R^1 \tilde\pi_*^G
\Omega_{\tilde C/k}^\vee) \to \H^0\left(C/G, \left ( \R^1 \pi_*^G
\right )\hom(\Omega_{C/k}, \left( \varphi_* \O_{\tilde C} \right ) /
\O_C)\right). 
\end{multline}
Il s'agit maintenant d'étudier le morphisme 
$$\H^0(\tilde C/G, \R^1 \tilde\pi_*^G \Omega_{\tilde C/k}^\vee) \to \H^0\left(C/G, \left ( \R^1 \pi_*^G \right )\hom(\Omega_{C/k}, \left( \varphi_* \O_{\tilde C} \right ) / \O_C)\right).$$
Ce morphisme peut se décomposer, à l'aide de
\cite{CohomologyOfGroup}, 27.1, en une somme de termes <<locaux>> de la
forme 
\begin{multline*}
\bigoplus_{\beta\in\pi_0(\tilde C)/G} \H^0(C_\beta/G_\beta, \R^1 \tilde\pi_{\beta, *}^{D_\beta} \Omega_{C_\beta/k}^\vee) \to  \\ \bigoplus_{\beta\in\pi_0(\tilde C)/G} \H^0\left(C_\beta/G_\beta, \R^1 \pi_{\beta, *}^{D_\beta} \left(\varphi_\beta^* \hom(\Omega_{C/k}, \varphi_* \O_{\tilde C}/\O_C) \right)\right)
\end{multline*}
car 
$$\hom(\Omega_{C/k}, \varphi_* \O_{\tilde C}/\O_C) = \varphi_* \varphi^* \hom(\Omega_{C/k}, \varphi_* \O_{\tilde C}/\O_C)$$
(ce qui peut se voir à l'aide de la description locale explicite
donnée dans la section \ref{LienNormaliseeLocal}). Il suffit donc
d'étudier chacun de ces termes locaux. 
Soit $\beta \in \pi_0(\tilde C)/G$. La suite spectrale des extensions
de groupes $\textrm{III}_{2}^{p, q}:=\H^p(G_\beta, \H^q(T_\beta,M))$ (pour
tout $D_{\beta}$-module $M$) nous donne (comme $T_\beta$ agit trivialement et que
$G_\beta=D_\beta/T_\beta$) un diagramme commutatif à colonnes exactes 
\begin{small}
$$\xymatrix{
0\ar[d] & 0 \ar[d]\\
\H^0(C_\beta/G_\beta, \R^1 \tilde\pi_*^{G_\beta} \Omega_{C_\beta/k}^\vee)\ar[d]\ar[r]^(0.37){a_\beta} & 
\H^0\left(C_\beta/G_\beta, \R^1 \pi_{\beta, *}^{G_\beta} \left(\varphi_\beta^* \hom(\Omega_{C/k}, \varphi_* \O_{\tilde C}/\O_C) \right)\right) \ar[d]\\
\H^0(C_\beta/G_\beta, \R^1 \tilde\pi_*^{D_\beta} \Omega_{C_\beta/k}^\vee)\ar[d]\ar[r]^(0.37){b_\beta} & 
\H^0\left(C_\beta/G_\beta, \R^1 \pi_{\beta, *}^{D_\beta} \left(\varphi_\beta^* \hom(\Omega_{C/k}, \varphi_* \O_{\tilde C}/\O_C) \right)\right)\ar[d] \\
\H^0(C_\beta, \H^1(T_\beta, \Omega_{C_\beta/k}^\vee))\ar[r]^(0.37){c_\beta} &
\H^0\left(C_\beta, \H^1\left(T_\beta, \varphi_\beta^*
\hom(\Omega_{C/k}, \varphi_* \O_{\tilde C}/\O_C)\right)\right). 
}$$
\end{small}
De plus, on voit que $c_\beta$ est canoniquement isomorphe au morphisme
$$\H^1(T_\beta, \H^0(C_\beta, \Omega_{C_\beta/k}^\vee)) \to
\H^1(T_\beta, \H^0(C_\beta, \varphi_\beta^* \hom(\Omega_{\tilde C/k},
\varphi_* \O_{\tilde C}/\O_C)))$$ 
car $T_\beta$ agit trivialement sur $C_\beta$.

D'autre part, comme $C \to \Spec k$ est stable (en particulier,
$\Ext^0(\Omega_{C/k}, \O_C)=0$ d'après \cite{Deligne_Mumford} Lemma
1.4) la suite exacte \eqref{SuiteLienNormaliseGlobal} montre (en
passant aux sections globales) que le morphisme 
$$\bigoplus_{C_0 \in \pi_0(\tilde C)} \Ext^0(\Omega_{C_0}, \O_{C_0})
\to \H^0(C, \hom(\Omega_{C/k}, \left( \varphi_* \O_{\tilde C} \right )
/ \O_C))$$ est injectif. Par suite, le morphisme $c_\beta$ est
injectif. 
Le noyau de $b_\beta$ est donc égal au noyau de $a_\beta$ et la suite
exacte \eqref{SuiteNormalisee2} donne une suite exacte 
\begin{multline*}
\hskip -0.4cm 0 \to \H^0(C/G, \Fg) \to \H^0(C, \hom(\Omega_{C/k}, \left( \varphi_* \O_{\tilde C} \right ) / \O_C))^G \to  \H^0(C/G, \R^1 \pi_*^G \Omega_{C/k}^\vee) \\  \to \bigoplus_{\beta\in\pi_0(\tilde C)/G}\H^0(C_\beta/G_\beta, \R^1 \tilde\pi_*^{G_\beta} \Omega_{C_\beta/k}^\vee) \to \bigoplus_{\beta\in\pi_0(\tilde C)/G} \Im a_\beta \to 0.
\end{multline*}
Finalement, le faisceau $\R^1 \tilde\pi_*^{G_\beta}
\Omega_{C_\beta/k}^\vee$ étant à support discret (car $G_\beta$ agit
librement sur un ouvert dense de $C_\beta$), l'image de $a_\beta$ est
égale à $\H^0(C_\beta/G_\beta,
\Gg_\beta)$ (par définition de $\Gg_\beta$). On a donc la suite exacte
\eqref{SuiteNormaliseeTermeLocal}. \qed 

\begin{remark}
\label{LienNormaliseeGlobalGeneral}
On peut remarquer que dans le cas général des courbes à singularités
doubles ordinaires, donc $C\to \Spec k$ pas forcément stable ni même
propre, les deux suites \eqref{SuiteExacteFGlobale} et
\eqref{SuiteNormalisee2} sont encore exactes. 
\end{remark}

\begin{corollary}
\label{CalculExplicite2Termes}
Reprenons les notations précédentes. On a
\begin{multline*}
\dim_k \H^1(C/G, \ext^0_G(\Omega_{C/k}, \O_C))+\dim_k \H^0(C/G, \R^1
\pi_*^G \Omega_{C/k}^\vee)= \\ 
\dim_k \H^0(C,  \hom(\Omega_{C/k}, \left( \varphi_* \O_{\tilde C}
\right ) / \O_C))^G-\sum_{\beta \in \pi_0(\tilde C)/G}\dim_k
\Hom(\Omega_{C_\beta/k}, \O_{C_\beta})^{G_\beta}\\ 
\sum_{\beta \in \pi_0(\tilde
C)/G}\dim_k\Ext^1_{G_\beta}(\Omega_{C_\beta}, \O_{C_\beta}) -\dim_k
\H^0(C_\beta/G_\beta, \Gg_\beta) 
\end{multline*}
\end{corollary}

\preuve Pour tout $\beta \in \pi_0(\tilde C)/G$ on a un isomorphisme canonique
$$\ext^1_{G_\beta}(\Omega_{C_\beta/k}, \O_{C_\beta}) \cong \R^1
\pi_{\beta, *}^{G_\beta} \Omega_{C_\beta/k}^\vee$$ 
(ceci provient de la suite spectrale
$\textrm{II}_{2}^{p,q}:=\R^p\pi_{\beta, *}^{G_{\beta}}\ext^{q}
(\Omega_{C_{\beta}/k}, \O_{C_\beta}))$ et de la lissité de $C_{\beta}$).
La suite spectrale $\textrm{I}_{2}^{p, q}:=\H^p (C_{\beta}/G_{\beta},
\ext^{q} (\Omega_{C_{\beta}/k}, \O_{C_\beta}))$ fournit donc une suite exacte
\begin{multline*}
0 \to \H^1(C_\beta, \ext^0_{G_\beta}(\Omega_{C_\beta/k},
\O_{C_\beta})) \to \Ext^1_{G_\beta}(\Omega_{C/k}, \O_C) \to \\
\H^0(C_\beta/G_\beta, \R^1 \pi_{\beta, *}^{G_\beta}
\Omega_{C_\beta/k}^\vee) \to 0. 
\end{multline*}
On peut alors décomposer ces groupes en prenant un représentant dans
chaque orbites sous l'action de $G$ via \cite{CohomologyOfGroup},
27.1, ce qui nous donne un isomorphisme
$$\Ext^0_G(\Omega_{\tilde C/k}, \O_{\tilde C}) \cong \bigoplus_{\beta
\in \pi_0(\tilde C)/G} \Ext^0(\Omega_{C_\beta/k},
\O_{C_\beta})^{G_\beta}$$ 
puis
$$\H^1(\tilde C/G, \ext^0_G(\Omega_{\tilde C/k}, \O_{\tilde C})) \cong
\bigoplus_{\beta \in \pi_0(\tilde C)/G}\H^1(C_\beta/G_\beta,
\ext^0_{G_\beta}(\Omega_{C_\beta/k}, \O_{C_\beta})).$$ 
Les suites exactes \eqref{SuiteNormaliseeTermeGlobal} et 
\eqref{SuiteNormaliseeTermeLocal} permettent alors de conclure. \qed 

\begin{exemple}
\label{ExempleGlobalCyclique}
Soient $C\to \Spec k$ une courbe stable sur un corps algébriquement
clos de caractéristique $p$ et $G$ un $p$-groupe cyclique agissant
fidèlement sur $C$. On ne suppose pas que le groupe $G$ agit librement
sur un ouvert dense de $C$. Supposons que les composantes connexes de
$\tilde C$ sont de genre $\ge 2$ et qu'en chaque point singulier $\p
\in C$, le groupe $D_\p$ agit sans permuter les branches. Pour tout
$\alpha \in C_{sing}/G$, choisissons un point $\p_\alpha$ dans la
classe $\alpha$ et notons $(\dif_\alpha, \dif'_\alpha)$ la différente
de $D_{\p_\alpha}$. Notons $\tilde C_{ram}$ l'ensemble des points $\q$
de $\tilde C$ tels que $G_\q$ (i.e. l'image du stabilisateur $D_\q$
dans $\Aut(\Spec \O_{\tilde C, \q})$) soit non trivial, choisissons
pour tout $\gamma \in \tilde C_{ram}/G$ un élément $\q_\gamma$ dans la
classe $\gamma$ et notons $\dif_{\gamma}$ la différente de
$G_{\q_\gamma}$ (et non $D_{\q_\gamma}$) si $G_{\q_\gamma} \not =
\{Id\}$ et $0$ sinon. Le corollaire \ref{CalculExplicite2Termes}, le
théorème \ref{CalculGlobalCasLisse} et l'exemple
\ref{CalculExpliciteG} montrent qu'alors 
\begin{multline*}
\dim_k \H^1(C/G, \ext^0_G(\Omega_{C/k}, \O_C))+\dim_k \H^0(C/G, \R^1
\pi_*^G \Omega_{C/k}^\vee)=2 \left|C_{sing}/G\right|+ \\\sum_{\beta
\in \pi_0(\tilde C/G)} 3p_a(C_\beta/G)-3+\sum_{\gamma \in \tilde
C_{ram}/G} \left\lfloor{ 2\dif_{\gamma} \over |D_{\q_\gamma}|} \right
\rfloor-\\\left|\left\{\substack{\alpha \in C_{sing}/G \\
2\dif_\alpha+1 \not = 0 \mod |G_{\p_\alpha}|\\ \textrm{ou }
G_{\p_\alpha} \not =
D_{\p_\alpha}}\right\}\right|-\left|\left\{\substack{\alpha \in
C_{sing}/G \\ 2\dif'_\alpha+1 \not = 0 \mod |G_{\p_\alpha}|\\
\textrm{ou } G_{\p_\alpha} \not = D_{\p_\alpha}}\right\}\right| 
\end{multline*}
le calcul de $\dim_k \H^0(C,  \hom(\Omega_{C/k}, \left( \varphi_*
\O_{\tilde C} \right ) / \O_C))^G$ étant immédiat (cf. description
locale dans la section \ref{LienNormaliseeLocal}). 
\end{exemple}

\subsubsection{Épaississement des points doubles}

Soient $k$ un corps algébriquement clos, $C \to \Spec k$ une courbe
stable et $G$ un groupe agissant sur $C$. 
Le but de cette section est de préciser l'image du morphisme
$$\Ext^1_G(\Omega_{C/k}, \O_C) \overset{{\textit{ép}}}{\to} \H^0(C,
\ext^1(\Omega_{C/k}, \O_C))^G.$$ 
Pour tout $\alpha \in C_{sing}/G$, choisissons un représentant
$\p_\alpha$ dans la classe $\alpha$ et notons $D_\alpha$ son
stabilisateur, $T_\alpha$ le noyau de l'homomorphisme $D_\alpha\to
\Aut(\Frac(\O_{C, \p_\alpha}))$ et
$G_\alpha:=D_\alpha/T_\alpha$. D'après \cite{CohomologyOfGroup} 27.1,
on a une décomposition
$$\H^0(C, \ext^1(\Omega_{C/k}, \O_C))^G=\bigoplus_{\alpha \in
C_{sing}/G} \Ext^1(\Omega_{C/k, \p_\alpha}, \O_{C,
\p_\alpha})^{D_\alpha}.$$ 

Le théorème \ref{RelevementActionSingulier} suggère de répartir les
points singuliers selon la nature des relèvements. Cela est fait dans
la définition suivante. 

\begin{definition}
Soient $\p$ un point singulier de $C$ et $D_\p$ son stabilisateur. Le
point $\p$ sera dit \textit{inconditionnellement relevable}
(resp. \textit{conditionnellement relevable}, resp. \textit{non
relevable}) si c'est le cas de l'action de $D_\p$ sur $\Spec \hat
\O_{C, \p}$ (cf. définition \ref{ActionRelevable}). 
\end{definition}

Grâce au théorème \ref{RelevementActionSingulier}, on voit que le
conducteur en un point conditionnellement relevable est de la forme
$(1, m)$ ou $(m, 1)$. Les points ayant un conducteur de cette forme
vont jouer un rôle particulier.  

Notons $\varphi:\tilde C \to C$ le morphisme de normalisation.
Définissons le sous-ensem\-ble $C_{cond \ 1}$ de
$\varphi^{-1}(C_{sing})$ de la manière suivante : 
$$C_{cond \ 1} :=\left\{
 \begin{array}{l}
\q_{1} \in \varphi^{-1}(C_{sing})|\varphi(\q_{1})\textrm{ est relevable et si } \q_{2} \textrm{ est l'autre point de} \\ \varphi^{-1}(\varphi(\q_{1}))
\textrm{ alors le conducteur de } T_{\q_1} \textrm{ sur } \Spec \O_{\tilde C, \q_{2}} \textrm { est } 1
\end{array}\right\}$$
(en particulier, les points de $\varphi (C_{cond\ 1})$ sont
conditionnellement relevable d'après le lemme \ref{ActionSurAnneauDeFraction}).
On définit alors deux diviseurs effectifs sur $\tilde C$  :

\begin{center}
$D_{sing}:=\displaystyle{\sum_{\q \in \pi^{-1}(C_{sing})}} \q,$ \ 
$D_{cond \ 1}:=\displaystyle{\sum_{\q_{1} \in C_{cond \ 1}}} \q_{1}. $ \
\end{center}
 
Pour tout $\beta \in \pi_0(\tilde C)/G$, choisissons un représentant
$C_\beta$ dans la classe de $\beta$ et notons $D_\beta$ son groupe
d'inertie, $T_\beta$ le noyau de l'homomorphisme $D_\beta \to
\Aut(\C_\beta)$ et $G_\beta:=D_\beta/T_\beta$.  
Notons $$\Na_\beta:=\H^1(T_\beta, \H^0(C_\beta,
\Omega_{C_\beta/k}^\vee(-D_{sing}+D_{cond \ 1})))^{G_\beta}.$$ 
Si on suppose que $\deg(D_{sing}-D_{cond \ 1}) > 2-2p_a(C_\beta)$ (ce
qui est le cas général et est toujours vérifié si $p_a(C_\beta) \ge
2$) on voit que $\Na_\beta=0$.  

\begin{theorem}
\label{EpaississementPointdouble}
Soient $C\to \Spec k$ une courbe stable sur un corps algébriquement
clos de caractéristique $p$ et $G$ un $p$-groupe fini agissant
fidèlement sur $C$. Pour tout $\alpha \in C_{sing}/G$ choisissons un
représentant $\p_\alpha$ dans la classe de $\alpha$. Reprenons les
notations ci-dessus. On a une suite exacte canonique   
$$
0 \to \bigoplus_{\substack{\alpha \in C_{sing}/G \\ \p_\alpha \textrm{
inconditionnellement}\\\textrm{relevable}}} \Ext^1(\Omega_{C,
\p_\alpha}, \O_{C, \p_\alpha}) \to \Im {\textit{ép}} \overset{a}{\to}
\bigoplus_{\beta \in \pi_0(\tilde C)/G} \Na_\beta. 
$$
\end{theorem}

Contrairement aux démonstrations précédentes, la démonstration de ce
thé\-orème n'est pas purement cohomologique mais repose également sur
des techniques de recollement (cf. \cite{Papier1}, Proposition 4.7).

\preuve Dans un premier temps, nous allons regarder le cas des points
inconditionnellement relevables. Ensuite nous construirons le
morphisme $a$ puis nous montrerons que la suite est exacte. 

\medskip

\noindent{\it Première étape : les points inconditionnellement relevables.}

Soit $\alpha \in C_{sing}/G$. Supposons que $\p_\alpha$ soit
inconditionnellement relevable.  
Nous allons construire une déformation $G$-équivariante $\C$ de $C$
qui est non topologiquement triviale en les points de $G.\p_\alpha$ et
topologiquement triviale en les points de $C_{sing}\setminus
G.\p_\alpha$. Ceci suffira à montrer qu'on a une inclusion
$$\bigoplus_{\substack{\alpha \in C_{sing}/G \\ \p_\alpha \textrm{
inconditionnellement}\\\textrm{relevable}}} \Ext^1(\Omega_{C,
\p_\alpha}, \O_{C, \p_\alpha}) \subset \Im {\textit{ép}}.$$

Notons $X=\Spec \hat\O_{C, \p_\alpha}$, $\eta$ et $\xi$ les deux
points génériques de $X$, $G_\eta$ et $G_\xi$ les groupes
d'automorphismes induits par $D_{\p_\alpha}$ sur $\kappa(\eta)$ et
$\kappa(\xi)$. Posons $C_\xi=\overline{\{\xi\}}$,
$C_\eta=\overline{\{\eta\}}$ et $U=\left(C_\xi \cup C_\eta\right)
\setminus C_{sing}$ (les adhérences étant prises dans $C$ en
identifiant les points $\eta$ et $\xi$ à leur image dans $C$). En
particulier, on a une identification naturelle de $G_\eta$
(resp. $G_\xi$) avec un sous-groupe de $\Aut ( C_\eta)$
(resp. $\Aut (C_\xi)$). 

Par définition d'un point inconditionnellement relevable, il existe
une déformation $D_{\p_\alpha}$-équivariante $\X$ de $X$ telle que $G_\eta$
et $G_\xi$ soient les groupes d'auto\-morphimes induits par
$D_{\p_\alpha}$ sur $\Spec \O_{\X, \eta}$ et $\Spec \O_{\X, \xi}$ (on
identifie les points génériques de $\X$ à ceux de $X$ via l'immersion
$X \to \X$). Comme il existe une unique déformation $G_\xi$
(resp. $G_\eta$) équivariante de $\Spec\kappa(\xi)$
(resp. $\Spec\kappa(\eta)$) (car les morphismes $\Spec
\kappa(\xi) \to \Spec \kappa(\xi)^{G_\xi}$ et $\Spec \kappa(\eta) \to
\Spec \kappa(\eta)^{G_\eta}$ sont étales), il existe un isomorphismes
$G_\eta$-équivariant  
$$\varphi_\eta:\Frac(\hat\O_{C_\eta, \p_\alpha}) \otimes_k k[\epsilon] \to \O_{\X, \eta}$$
et un isomorphisme $G_\xi$-équivariant 
$$\varphi_\xi:\Frac(\hat\O_{C_\xi, \p_\alpha}) \otimes_k k[\epsilon] \to \O_{\X, \xi}.$$
En transportant ces isomorphismes à tous les points de $G.\p_\alpha$
via l'action de $G$ on peut construire, grâce à la proposition 4.7 de
\cite{Papier1}, une déformation équivariante $(\C, G)$ de $(C, G)$ qui
induit la déformation triviale de $C\setminus G.\p_\alpha$ et une
déformation non topologiquement triviale de $\Spec \widehat\O_{C,
\p_\alpha}$. 

En conclusion, on a une injection
$$\bigoplus_{\substack{\alpha \in C_{sing}/G \\ \p_\alpha \textrm{
inconditionnellement relevable}}} \Ext^1(\Omega_{C, \p_\alpha}, \O_{C,
\p_\alpha}) \to \Im {\textit{ép}}.$$ 

\medskip

\noindent{\it Deuxième étape : construction de $a$.}

On peut remarquer que l'image de ${\textit{ép}}$ s'identifie avec l'image de
$$\H^0(C/G, \ext^1_G(\Omega_{C/k}, \O_C)) \to \H^0(C,
\ext^1(\Omega_{C/k}, \O_C))^G.$$ C'est ce dernier morphisme que nous
allons étudier. 
Reprenons les notations de la section \ref{LienNormalisee} concernant
$C_\beta$, $\pi_\beta$, $\varphi_\beta$, ... (cf. diagramme \eqref{DefinitionMorphismesGlobaux})

La suite spectrale $\textrm{II}_{2}^{p, q}:=\R^p \pi_{*}^{G}\ext^{q}
(\Omega_{C/k}, \O_C)$ fournit une suite exacte
$$0 \to \R^1 \pi_*^G \Omega_{C/k}^\vee \to \ext^1_G(\Omega_{C/k},
\O_C) \to \pi_*^G \ext^1(\Omega_{C/k}, \O_C).$$ 
En particulier, si $U$ désigne l'image dans $C/G$ du lieu lisse de $C$, on a un isomorphisme 
$$\left(\R^1 \pi_*^G \Omega_{C/k}^\vee \right) |_U \cong \ext^1_G(\Omega_{C/k}, \O_C) |_U.$$
En décomposant $\pi^{-1} (U)$ selon les orbites sous l'action de $G$ on peut
écrire, grâce à \cite{CohomologyOfGroup} 27.1, un isomorphisme 
$$\left(\R^1 \pi_*^G \Omega_{C/k}^\vee \right) |_U\cong\bigoplus_{\beta \in
\pi_0(\tilde C)/G}\left( \tilde \varphi_{\beta, *} \R^1\pi_{\beta,
*}^{D_\beta} \Omega_{C_\beta/k}^\vee \right) |_U.$$ 
Notons $\theta$ le morphisme composé
\begin{multline*}
\ext^1_G(\Omega_{C/k}, \O_C) \to \ext^1_G(\Omega_{C/k}, \O_C)|_U
\overset{\sim}{\to} \left(\R^1 \pi_*^G \Omega_{C/k}^\vee \right)
|_U\overset{\sim}{\to} \\ \bigoplus_{\beta \in \pi_0(\tilde C)/G}
\left( \tilde \varphi_{\beta, *} \R^1\pi_{\beta, *}^{D_\beta}
\Omega_{C_\beta/k}^\vee\right) |_U 
\end{multline*}
et $\theta'$ le composé de $\theta$ avec l'homomorphisme de changement de groupes
$$\bigoplus_{\beta \in \pi_0(\tilde C)/G} \left(\tilde \varphi_{\beta, *}
\R^1\pi_{\beta, *}^{D_\beta} \Omega_{C_\beta/k}^\vee \right) |_U \to
\bigoplus_{\beta \in \pi_0(\tilde C)/G} \left( \tilde \varphi_{\beta, *}
\pi_{\beta, *}^{G_\beta} \H^1(T_\beta, \Omega_{C_\beta/k}^\vee)\right) |_U.$$ 

Le morphisme $\theta$ peut se décrire, localement, de manière explicite. En effet,
choisissons un point $\p\in C_{sing}$ et fixons une composante
$C_\beta \in \pi_0(\tilde C)$ telle que $\varphi^{-1}(\p)\cap C_\beta
\not = \emptyset$. Soit $\q \in \varphi^{-1}(\p)\cap C_\beta$. Toute
déformation $G$-équivariante du premier ordre $(\C\to\Spec
k[\epsilon], G)$ de $C$ induit une déformation $D_\p$-équivariante de
$\Spec \hat \O_{\C, \p}$. Pour tout $\sigma \in D_\p$, notons
$\sigma_\epsilon$ l'automorphisme de $\Spec \hat \O_{\C, \p}$
induit. Le lemme \ref{RelevementMorphisme} (dont on reprend les
notations, l'uniformisante $x$ correspondant à une uniformisante de
$C_\beta$ en $\q$) fournit alors la forme de $\sigma_\epsilon$ pour
tout $\sigma$. L'image par $\theta$ de $\C$ est alors donnée
localement (via le lemme \ref{ActionSurAnneauDeFraction}) par le
morphisme de cochaînes $\phi$ défini pour tout $\sigma \in D_\p\cap
D_\beta$ par  
\begin{equation}
\label{psiexplicite}
\phi:\sigma \mapsto {1 \over \der{x} \sigma(x)}\left( \lambda
v_{\sigma, 1}+ h_{\sigma, 0}(x) \right )\d x. 
\end{equation}
En particulier, on voit que pour tout $\sigma\in D_\p \cap D_\beta$ la série de Laurent
$\phi(\sigma)(dx)$
n'a pas de pôle. Le morphisme $\theta$ se factorise donc par  
$$\theta_1:\ext^1_G(\Omega_{C/k}, \O_C) \to \bigoplus_{\beta \in
\pi_0(\tilde C)/G} \tilde \varphi_{\beta, *} \R^1\pi_{\beta,
*}^{D_\beta} \Omega_{C_\beta/k}^\vee.$$ De même, on voit que le
morphisme $\theta'$ se factorise par  
$$\theta'_1:\ext^1_G(\Omega_{C/k}, \O_C)\to \bigoplus_{\beta \in \pi_0(\tilde C)/G}
\tilde \varphi_{\beta, *}\pi_{\beta, *}^{G_\beta} \H^1(T_\beta,
\Omega_{C_\beta/k}^\vee)$$ et l'image par $\theta'_1$ de $\R^1 \pi_*^G
\Omega_{C/k}^\vee$ (dont les éléments correspondent aux déformations
topologiquement triviales, i.e. $\lambda = 0$) est dans  
$$\bigoplus_{\beta \in \pi_0(\tilde C)/G}
\tilde \varphi_{\beta, *}\pi_{\beta, *}^{G_\beta} \H^1(T_\beta,
\Omega_{C_\beta/k}^\vee(-D_{sing})).$$ 
De plus, si $\sigma \in T_\beta$ et que  $\lambda$ est non nulle
(c'est-à-dire que la déformation équivariante de $\Spec \hat \O_{C,
\p}$ induite par $\C$ n'est pas topologiquement triviale) alors
$\phi(\sigma)(dx)$ possède un zéro si et seulement si $\p \not\in
C_{cond \ 1}$. 
Par suite, on voit que le morphisme $\theta'_1$ se factorise par 
$$\theta'_2:\ext^1_G(\Omega_{C/k}, \O_C)\to\bigoplus_{\beta \in \pi_0(\tilde C)/G}
\tilde \varphi_{\beta, *}\pi_{\beta, *}^{G_\beta} \H^1(T_\beta,
\Omega_{C_\beta/k}^\vee(-D_{sing}+D_{cond \ 1})).$$ 
En résumé, on a un diagramme commutatif
$$\xymatrix{
\R^1 \pi_*^G \Omega_{C/k}^\vee \ar@{^(->}[d] \ar[r]^(.32){\theta'_2} &
\bigoplus_{\beta \in \pi_0(\tilde C)/G} \tilde \varphi_{\beta,
*}\pi_{\beta, *}^{G_\beta} \H^1(T_\beta,
\Omega_{C_\beta/k}^\vee(-D_{sing})) \ar@{^(->}[d] \\ 
\ext^1_G(\Omega_{C/k}, \O_C) \ar[r]_(.27){\theta'_2} & \bigoplus_{\beta \in \pi_0(\tilde C)/G}
\tilde \varphi_{\beta, *}\pi_{\beta, *}^{G_\beta} \H^1(T_\beta,
\Omega_{C_\beta/k}^\vee(-D_{sing}+D_{cond \ 1}))
}$$
En considérant la suite exacte longue de cohomologie (des faisceaux)
associée à chacune des colonnes, on trouve un diagramme commutatif à
colonnes exactes 
\begin{footnotesize}
$$\xymatrix{
0 \ar[d] & 0\ar[d]\\
\H^0(C/G, \R^1 \pi_*^G \Omega_{C/k}^\vee) \ar[d]
\ar[r]^(.37){\H^0(\theta'_2)} & \displaystyle\bigoplus_{\beta \in
\pi_0(\tilde C)/G} \H^1(T_\beta, \H^0(C_\beta,
\Omega_{C_\beta/k}^\vee(-D_{sing})))^{G_\beta} \ar[d] \\ 
\H^0(C/G, \ext^1_G(\Omega_{C/k}, \O_C))
\ar[d]\ar[r]^(.37){\H^0(\theta'_2)} & \displaystyle\bigoplus_{\beta
\in \pi_0(\tilde C)/G} 
 \H^1(T_\beta, \H^0(C_\beta, \Omega_{C_\beta/k}^\vee(-D_{sing}+D_{cond
\ 1})))^{G_\beta} \\ 
\Im {\textit{ép}} \ar[d]&  \\
0 & 
}$$
\end{footnotesize}
Comme $C$ est une courbe stable, on voit que pour tout $\beta \in
\pi_0(\tilde C)$ on a $$\deg \Omega_{C_\beta/k}^\vee(-D_{sing}) < 0$$
donc $\H^0(C_\beta, \Omega_{C_\beta/k}^\vee(-D_{sing}))=0.$
On a donc un morphisme canonique
$$a:\Im {\textit{ép}} \to \bigoplus_{\beta \in \pi_0(\tilde C)/G}
  \H^1(T_\beta, \H^0(C_\beta, \Omega_{C_\beta/k}^\vee(-D_{sing}+D_{cond \ 1})))^{G_\beta}$$
induit par $\H^0(\theta'_2)$.
En résumé, on a un diagramme commutatif
\begin{footnotesize}
$$\xymatrix{
\H^0(C/G, \ext^1_G(\Omega_{C/k}, \O_C)) \ar@{->>}[d]
\ar[rdd]^{\H^0(\theta'_1)} \ar[r]^{\H^0(\theta_1)}
\ar@/_1cm/[dd]_{\H^0(\theta_2')} & \bigoplus_{\beta}
\H^0(C_\beta/G_\beta, \R^1 \pi_{\beta, *}^{D_\beta}
\Omega_{C_\beta/k}^\vee) \ar[dd]\\ 
\Im {\textit{ép}} \ar[d]_{a}\\
 \bigoplus_{\beta}
  \H^1(T_\beta, \H^0(C_\beta,
\Omega_{C_\beta/k}^\vee(-D_{sing}+D_{cond \ 1})))^{G_\beta}
\ar@{^(->}[r] & \bigoplus_{\beta} 
  \H^1(T_\beta, \H^0(C_\beta, \Omega_{C_\beta/k}^\vee))^{G_\beta}
}$$
\end{footnotesize}

\medskip

\noindent{\it Troisième étape : exactitude de la suite.} 

Soit $$f \in \bigoplus_{\substack{\alpha \in C_{sing}/G \\ \p_\alpha
\textrm{ inconditionnellement}\\\textrm{relevable}}} \Ext^1(\Omega_{C,
\p_\alpha}, \O_{C, \p_\alpha}).$$  
D'après la construction explicite de la première étape de la preuve,
on peut trouver une déformation équivariante $(\C, G)$ de $(C, G)$
telle que pour tout $\beta \in \pi_0(\tilde C)/G$ et tout ouvert $U
\subset C_\beta \cap \varphi_{\beta}^{-1} (C_{lisse})$, le groupe
$T_\beta$ agit trivialement sur la déformation équivariante de $(U,
D_\beta)$ induite par $(\C, G)$. Ainsi, l'image de $(\C, G)$ (et donc
de $f$) dans $\Na_\beta$ est nulle. Le morphisme $a$ se factorise donc
par le conoyau $\mathfrak{I}$ du morphisme 
$$\bigoplus_{\substack{\alpha \in C_{sing}/G \\ \p_\alpha \textrm{
inconditionnellement}\\\textrm{relevable}}} \Ext^1(\Omega_{C,
\p_\alpha}, \O_{C, \p_\alpha}) \to \Im {\textit{ép}}.$$ 
Montrons maintenant que le morphisme induit par $a$
$$\mathfrak{I} \to \bigoplus_{\beta \in \pi_0(\tilde C)/G} \H^1(T_\beta,
\H^0(C_\beta, \Omega_{C_\beta/k}^\vee(-D_{sing}+D_{cond \
1})))^{G_\beta}$$ 
est injectif.

Soit $f \in \mathfrak{I}$ un élément non nul. Alors par définition il existe une
déformation équivariante $(\C\to\Spec k[\epsilon], G)$ de $(C, G)$
dont l'image dans $\mathfrak{I}$ est $f$. Comme $f$ est non nul, il existe un
point $\p \in C_{sing}$ conditionnellement relevable tel que la
déformation $D_\p$ équivariante de $\Spec \widehat\O_{C, \p}$ induite
par $(\C, G)$ soit non topologiquement triviale. Le fait que le point
soit conditionnellement relevable impose que l'image de $f$ par
$a$ est non triviale (localement donc globalement). \qed 

\begin{remark}
Dans la cas où les stabilisateurs en chaque point singulier  sont
cycliques et vérifient certaines conditions (sur les conducteurs) il
est possibles de préciser l'image de l'application $a$ définie dans le
théorème précédent (cf. \cite{These} Chap II, Théorème 4.2.9).
\end{remark}

\begin{exemple}
Reprenons les hypothèses et notations de l'exemple
\ref{ExempleGlobalCyclique}. En particulier, comme on a supposé que
les composantes de $\tilde C$ sont de genre $\ge 2$, on a $\Na_\beta =
0$ pour tout $\beta \in \pi_0(\tilde C/G)$. Le théorème
\ref{EpaississementPointdouble} permet donc de voir que  
\begin{multline*}
\dim_k \Ext^1_G(\Omega_{C/k}, \O_C)=2 \left|C_{sing}/G\right|+
\\\sum_{\beta \in \pi_0(\tilde C/G)} 3p_a(C_\beta/G)-3+\sum_{\gamma
\in \tilde C_{ram}/G} \left\lfloor{ 2\dif_{\gamma} \over
|D_{\q_\gamma}|} \right \rfloor-\\\left|\left\{\substack{\alpha \in
C_{sing}/G \\ 2\dif_\alpha+1 \not = 0 \mod |G_{\p_\alpha}|\\
\textrm{ou } G_{\p_\alpha} \not =
D_{\p_\alpha}}\right\}\right|-\left|\left\{\substack{\alpha \in
C_{sing}/G \\ 2\dif'_\alpha+1 \not = 0 \mod |G_{\p_\alpha}|\\
\textrm{ou } G_{\p_\alpha} \not =
D_{\p_\alpha}}\right\}\right|+\left|\left\{\substack{\alpha \in
C_{sing}/G \\ \p_\alpha \textrm{
inconditionnellement}\\\textrm{relevable}}\right\}\right|. 
\end{multline*}

Dans le cas où l'action est libre sur un ouvert dense (en particulier,
il n'y a pas de points contitionnellement relevables) mais qu'on
ne suppose pas que les composantes irréductibles de $\tilde C$ sont de
genre $\ge 2$, alors la formule ci-dessus reste vraie.
\end{exemple}

\providecommand{\bysame}{\leavevmode\hbox to3em{\hrulefill}\thinspace}


\begin{thebibliography}{Mau03}

\bibitem[Bab69]{CohomologyOfGroup}
A.~Babakhanian, \emph{Cohomology of finite groups}, Queen's University,
  Kingston, Ont., 1969.

\bibitem[BM00]{BertinMezard}
J.~Bertin and A.~M{\'e}zard, \emph{D\'eformations formelles des rev\^etements
  sauvagement ramifi\'es de courbes alg\'ebriques}, Invent. Math.
  \textbf{141}:1 (2000), 195--238.

\bibitem[CK]{CornelissenKato}
G.~Cornelissen and F.~Kato, \emph{Equivariant deformation of Mumford curves and
  of ordinary curves in positive characteristic}, \mbox{arXiv:math.AG/0103207},
  à paraître dans Duke Math. Journal.

\bibitem[DM69]{Deligne_Mumford}
P.~Deligne and D.~Mumford, \emph{The irreducibility of the space of curves of
  given genus}, Inst. Hautes \'Etudes Sci. Publ. Math.:36 (1969), 75--109.

\bibitem[Gro57]{Tohoku}
A.~Grothendieck, \emph{Sur quelques points d'alg\`ebre homologique}, T\^ohoku
  Math. J. (2) \textbf{9} (1957), 119--221.

\bibitem[Ill71]{IllusieI}
L.~Illusie, \emph{Complexe cotangent et d\'eformations. {I}}, Springer-Verlag,
  Berlin, 1971, Lecture Notes in Mathematics, Vol. 239.

\bibitem[Ill72]{IllusieII}
L.~Illusie, \emph{Complexe cotangent et d\'eformations. {I}{I}},
  Springer-Verlag, Berlin, 1972, Lecture Notes in Mathematics, Vol. 283.

\bibitem[LL78]{Laudal_Lonsted}
O.~A. Laudal and K.~L{\o}nsted, \emph{Deformations of curves. {I}. {M}oduli for
  hyperelliptic curves}, Algebraic geometry (Proc. Sympos., Univ. Troms\o,
  Troms\o, 1977), Lecture Notes in Math., vol. 687, Springer, Berlin, 1978,
  pp.~150--167.

\bibitem[Mau03a]{Papier1}
S.~Maugeais, \emph{Rel\`evement des rev\^etements p-cycliques des courbes
  rationnelles semi-stables}, Math. Ann. \textbf{327}:2 (2003), 365--393.

\bibitem[Mau03b]{Papier3}
S.~Maugeais, \emph{Th\'eorie des d\'eformations \'equivariantes des morphismes
  localement d'intersections compl\`etes}, 2003,
  \mbox{http://arXiv.org/abs/math.AG/0310136}, preprint.

\bibitem[Mau03c]{These}
S.~Maugeais, \emph{D\'eformations \'equivariantes des courbes stables}, Th\`ese
  de doctorat, Universit\'e de Bordeaux I, 2003.

\bibitem[Sch68]{Deformation}
M.~Schlessinger, \emph{Functors of {A}rtin rings}, Trans. Amer. Math. Soc.
  \textbf{130} (1968), 208--222.

\bibitem[Ser62]{Corps_Locaux}
J.~P. Serre, \emph{Corps locaux}, Actualit\'es Sci. Indust., No. 1296. Hermann,
  Paris, 1962.

\bibitem[Tuf93]{Tuffery}
S.~Tuff{\'e}ry, \emph{D\'eformations de courbes avec action de groupe}, Forum
  Math. \textbf{5}:3 (1993), 243--259.

\bibitem[Wew]{DeformationWewers}
S.~Wewers, \emph{Formal deformation of curves with group scheme action},
  \mbox{arXiv:math.AG/0212145}, preprint.

\end{thebibliography}
\end{document}